\documentclass[11pt]{article}
\usepackage{mathrsfs}
\usepackage{amsmath}
\usepackage{amsthm}
\usepackage{amssymb}
\usepackage{amsfonts}
\usepackage{color,xcolor}
\usepackage{graphicx}
\usepackage{geometry}
\usepackage{manfnt}
\usepackage{amsmath}
\usepackage[pagewise]{lineno}
\usepackage[colorlinks,citecolor=blue,urlcolor=blue]{hyperref}
\usepackage{tikz}
\newtheorem{theorem}{Theorem}[section]

\newtheorem{lemma}[theorem]{Lemma}

\newtheorem{definition}[theorem]{Definition}
\newtheorem{remark}[theorem]{Remark}
\newtheorem{example}[theorem]{Example}

\begin{document}
\title{\bf Stochastic flows for SDEs with critically low regularity growing drifts in Lebesgue--H\"older--Dini spaces}

\author{Jinlong Wei$^1$,  Wei Wang$^2$, Guangying Lv$^3$\thanks{Corresponding author.} and Jinqiao Duan$^4$
\\ {\small \it $^1$School of Statistics and Mathematics, Zhongnan University of}\\
{\small \it Economics and Law, Wuhan 430073, P.R. China}
\\ {\small \tt  weijinlong.hust@gmail.com}
\\ {\small \it $^2$School of Mathematics, Nanjing University, Nanjing 210093, P.R. China}
\\ {\small \tt wangweinju@nju.edu.cn}
\\ {\small \it $^3$College of Mathematics and Statistics,
Nanjing University of Information}
\\{\small \it  Science and Technology, Nanjing 210044, P.R. China}\\ {\small \tt  gylvmaths@126.com}
\\  {\small \it $^4$College of Science, Great Bay University, Dongguan 523000, P.R. China}\\
 {\small \tt  duan@gbu.edu.cn}}

\date{}

\maketitle
\noindent{\hrulefill}
\vskip1mm\noindent{\bf Abstract} We study stochastic differential equations (SDEs) with critically low regularity drifts in a Lebesgue--H\"older--Dini class. Under suitable slowly varying and H\"older--Dini-type assumptions on the modulus of continuity, we prove that SDEs with linearly growing drifts in the critical class
$L^p\big([0,T];{\mathcal C}^{\frac{2}{p}-1,\rho}_{l,d}(\mathbb{R}^n;\mathbb{R}^n)\big)$, $p\in(1,2]$,
generate a unique stochastic flow of homeomorphisms. The proof relies on a maximal regularity theory for parabolic equations in non-divergence form with Lebesgue--H\"older--Dini coefficients, which enables us to apply the It\^o--Tanaka trick and reduce the original SDE to one with Lipschitz coefficients.


\vskip0mm\noindent{\hrulefill}

\section{Introduction}\label{sec1}\setcounter{equation}{0}
Strong well-posedness for stochastic differential equations (SDEs) with low regularity coefficients
constitutes a central topic in stochastic analysis. In the classical It\^o theory, this issue is well understood under Lipschitz assumptions on the coefficients. In critical low regularity regimes, however, especially when the drift belongs only to borderline spaces, the problem becomes substantially more delicate.

\smallskip
In this paper we study the stochastic differential equation in ${\mathbb R}^n$
\begin{align}\label{1.1}
dX_{s,t}(x)=b(t,X_{s,t}(x))dt+dW_t,\quad t\in(s,T],\ \ X_{s,s}=x\in{\mathbb R}^n,
\end{align}
where $s\in[0,T]$ and $\{W_t\}_{t\in[0,T]}$ is an $n$-dimensional standard Wiener process. Classical results of It\^o \cite{Ito} yield strong well-posedness for Lipschitz drifts, while Veretennikov \cite{Ver} established strong well-posedness for bounded measurable drifts. When the drift is merely integrable in time or has low spatial regularity, the problem becomes much subtler.

\smallskip
A milestone result of Krylov and R\"ockner \cite{KR} established strong well-posedness under the subcritical Ladyzhenskaya--Prodi--Serrin (LPS) condition
\[
b\in L^p([0,T];L^q({\mathbb R}^n;{\mathbb R}^n)),\quad p,q\in[2,\infty],\ \ \tfrac{2}{p}+\tfrac{n}{q}<1.
\]
This theory was later extended to non-constant diffusion coefficients
\cite{XX,ZY,Zha05,Zha11}. By contrast, the critical balance $2/p+n/q=1$ remains a major challenge. Partial progress was obtained by R\"ockner and Zhao \cite{RZ2} for $p,q>2$ and $n\geq 3$, and by Beck, Flandoli, Gubinelli and Maurelli \cite{BFGM} for $p=2$, $q=\infty$ with a diffuse initial datum; see also \cite{FT,KS,KM,Kin,Kry21-1,Kry21-2,Kry21-3,Kry23-1,Kry23-2,Nam,RZ1,WLWu,YZ}.
More recently, Krylov \cite{Kry25} proved conditional pathwise
uniqueness for equations with Morrey drift and VMO diffusion,
covering the case $p=2$, $q=\infty$. However, this result does not by itself establish strong
well-posedness under the sole assumption
$b\in L^2([0,T];L^\infty(\mathbb R^n;\mathbb R^n))$.

\smallskip
Motivated by the classical result in \cite{Ito}, as well as by the critical results in \cite{BFGM,Kry25},
we seek intermediate function spaces between $L^1([0,T];\mathrm{Lip}({\mathbb R}^n))$
and $L^2([0,T];L^\infty({\mathbb R}^n))$ in which \eqref{1.1} remains strongly well posed.
A natural candidate is $L^p([0,T];{\mathcal C}^\alpha_l({\mathbb R}^n))$ with $p\in(1,2]$ and $\alpha\in[0,1)$,
where ${\mathcal C}^\alpha_l$ denotes the space of locally H\"older continuous functions for $\alpha>0$,
and uniformly continuous functions for $\alpha=0$.
By scaling, the analogue of the LPS condition takes the form
\begin{align}\label{1.2}
b\in L^p([0,T];{\mathcal C}^\alpha_l({\mathbb R}^n;{\mathbb R}^n)),\quad
p\in(1,2],\ \ \alpha\in[0,1),\ \ \tfrac{2}{p}-\alpha\le1.
\end{align}
The strong well-posedness problem for \eqref{1.1} under \eqref{1.2} is therefore as natural as in the classical LPS framework, and comparably delicate at the critical threshold.

\smallskip
In the subcritical case $2/p-\alpha<1$, strong well-posedness has been established in
\cite{GG,TDW,WHY1}.
Subsequently, \cite{WHY2} showed that, at the critical threshold, well-posedness still holds if $b$ enjoys additional Lorentz regularity,
\[
b\in L^{p,1}\big([0,T];{\mathcal C}^{\frac{2}{p}-1}_l({\mathbb R}^n;{\mathbb R}^n)\big), \quad p\in (1,2).
\]
However, the critical case $2/p-\alpha=1$ remains open when $b$ is assumed to be only $L^p$-integrable in time.

\smallskip
To study this critical regime, we employ the It\^o--Tanaka trick, also known as the Zvonkin transform \cite{Zvo}.
This approach is based on solving a suitable parabolic equation with drift and transforming \eqref{1.1}
into an SDE with more regular coefficients.
At the critical threshold, this requires a maximal regularity theory in function spaces associated with H\"older--Dini moduli of continuity.

\smallskip
The purpose of this paper is to address this critical case by introducing a new functional framework based on H\"older--Dini moduli that are slowly varying at the origin, inspired by \cite{WZ}.
Within this framework, we establish strong well-posedness for \eqref{1.1} with drift in critical Lebesgue--H\"older--Dini classes, and further prove that the solution family generates a stochastic flow of homeomorphisms.
In particular, for $p=2$ we obtain a partial answer to the above open problem under a natural replacement of $L^\infty$ by a locally Dini space.

\subsection{Main result}\label{sec1.1}
We now state the main result. Since it is formulated in terms of stochastic flows, we first recall the following definition.
\begin{definition}[\cite{Kun90}, p.~114]\label{def1.1}
A stochastic flow of homeomorphisms on a stochastic basis
$(\Omega, \mathcal{F},{\mathbb P},$ $ \{{\mathcal F}_t\}_{t\in [0,T]})$ associated with SDE \eqref{1.1}
is a map
\[
(s,t,x,\omega)\mapsto X_{s,t}(x,\omega), \quad 0\le s\le t\le T,\; x\in\mathbb{R}^n,\; \omega\in\Omega,
\]
with values in $\mathbb R^n$, such that:
\begin{enumerate}
\item[(i)] For each $s\in [0,T]$ and $x\in \mathbb{R}^n$, the process
$\{X_{s,t}(x),t\in [s,T]\}$ is a continuous
$\{\mathcal{F}_{s,t}\}_{s\le t\le T}$-adapted solution of SDE \eqref{1.1}.
\item[(ii)] ${\mathbb P}$-a.s., $X_{s,t}(\cdot)$ is a homeomorphism for all $0\le s\le t\le T$, and both $X_{s,t}(x)$ and $X_{s,t}^{-1}(x)$ are continuous in $(s,t,x)$.
\item[(iii)] ${\mathbb P}$-a.s., the flow property holds: $X_{s,t}(x)=X_{r,t}(X_{s,r}(x))$ for all
$0\le s\le r \le t \le T$, $x\in {\mathbb R}^n$, and $X_{s,s}(x)=x$.
\end{enumerate}
Here $\mathcal F_{s,t}$ denotes the completed sigma-field generated by $\{W_r-W_s:s\le r\le t\}$.
\end{definition}

Our main result is given as follows.
\begin{theorem}\label{the1.2}
Let $b\in L^p([0,T];{\mathcal C}^{\frac{2}{p}-1,\rho}_{l,d}(\mathbb{R}^n;\mathbb{R}^n))$ with $p\in (1,2]$,
and assume that $\rho^{\frac{2p}{5p-2}}$ is a Dini function and that $\rho$ is slowly varying at zero.
Then SDE \eqref{1.1} generates a unique stochastic flow of homeomorphisms
$\{X_{s,t}(x),\; t\in [s,T]\}$. Moreover, for every $R>0$, all $x,x^\prime\in \mathbb{R}^n$ with $|x|,|x^\prime|\le R$, all $0\le s\le t\le T$, $0\le s^\prime\le t^\prime\le T$, and every $q\geq 2$,
\[
\mathbb{E}|X_{s,t}(x)-X_{s^\prime,t^\prime}(x^\prime)|^q
\le C\bigl(|x-x^\prime|^q+|s-s^\prime|^{(1-\frac{1}{p})q}+|t-t^\prime|^{(1-\frac{1}{p})q}\bigr).
\]
\end{theorem}

\begin{remark}\label{rem1.3}
Theorem \ref{the1.2} shows that a stochastic flow of homeomorphisms is still generated at the critical threshold $\alpha=2/p-1$, provided that the H\"older modulus is refined by a slowly varying Dini factor. In this way, the theorem extends the subcritical well-posedness theory to the borderline regime and yields a new functional framework for critical SDEs with low regularity drift. Moreover, the theorem establishes a quantitative $L^q$-estimate for the dependence of the flow on the spatial variable and on the initial and terminal times.
\end{remark}

For $1<p<2$ and $\alpha=2/p-1$, Wei, Hu and Yuan \cite{WHY2}
proved that SDE~\eqref{1.1} generates a unique stochastic flow
of homeomorphisms when the drift is $L^{p,1}$-integrable in
time and has spatial modulus $r^\alpha$.
Theorem~\ref{the1.2} obtains the same conclusion by relaxing
the time integrability to $L^p$ while requiring the finer
spatial modulus $r^\alpha\rho(r)$. The following example gives
a drift satisfying the assumptions of Theorem~\ref{the1.2}
whose time dependence lies outside $L^{p,1}$ and whose spatial
profile has no H\"older exponent greater than $\alpha$ at
the origin.
\begin{example}\label{exa1.4}
Let $n=T=1$, fix $1<p\le2$, and set $\alpha=2/p-1$.
The function $\rho(r)=(\log(e/r))^{-3}$, $0<r\le1$, is slowly
varying at zero, and $\rho^{2p/(5p-2)}$ satisfies the Dini condition
since
\[
  \int_0^1 \rho(r)^{\frac{2p}{5p-2}}\frac{dr}{r}
  =\int_1^\infty \tau^{-\frac{6p}{5p-2}}d\tau<\infty.
\]
Choose $r_0\in(0,1)$ sufficiently small that
$r\mapsto r^\alpha\rho(r)$ is concave on $(0,r_0]$.
Define $h:\mathbb R\to[0,\infty)$ by
\[
 h(x)=
 \begin{cases}
  0, & x=0,\\
  |x|^\alpha\rho(|x|), & 0<|x|\le r_0,\\
  r_0^\alpha\rho(r_0), & |x|>r_0.
 \end{cases}
\]
Let
\begin{equation}\label{1.3}
 b(t,x)=b_1(t)h(x),\qquad
 b_1(t)=\frac{t^{-1/p}}{\log(e^2/t)},
 \qquad (t,x)\in(0,1]\times\mathbb R,
\end{equation}
and set $b(0,x)=0$. Since $b_1$ is decreasing on $(0,1]$, it coincides with
its decreasing rearrangement. The substitution
$\tau=\log(e^2/t)$ gives
\[
 \int_0^1 b_1(t)^pdt
 =\int_2^\infty\tau^{-p}d\tau<\infty,
 \qquad
 \int_0^1 t^{\frac1p-1}b_1(t)dt
 =\int_2^\infty\frac{d\tau}{\tau}=\infty.
\]
Thus $b_1\in L^p(0,1)\setminus L^{p,1}(0,1)$. On the other hand, $h$ is bounded and even, and its
restriction to $[0,\infty)$ is nondecreasing and concave
with $h(0)=0$. This restriction is therefore subadditive.
Since $r\mapsto r^\alpha\rho(r)$ is increasing on $(0,1]$,
we also have
\[
 h(r)\le r^\alpha\rho(r),\qquad 0<r\le1.
\]
Consequently, for $x,y\in\mathbb R$ with $0<|x-y|\le1$,
\[
 |h(x)-h(y)|
 \le h\bigl(\big||x|-|y|\big|\bigr)
 \le h(|x-y|)
 \le |x-y|^\alpha\rho(|x-y|).
\]
Together with the boundedness of $h$, this yields $b\in L^p([0,1];
 \mathcal C^{\alpha,\rho}_{l,d}(\mathbb R))$.

\smallskip
By Theorem~\ref{the1.2}, SDE~\eqref{1.1} with drift~\eqref{1.3} generates a unique stochastic flow of homeomorphisms. However, for every $t\in(0,1]$ and every $\beta>\alpha$,
\[
  \frac{|b(t,r)-b(t,0)|}{r^\beta}
  =\frac{b_1(t)}{r^{\beta-\alpha}(\log(e/r))^3}
  \longrightarrow\infty
  \qquad{\rm as}\quad r\downarrow0.
\]
Thus $b(t,\cdot)$ fails to satisfy any subcritical H\"older
condition at the given $p$. In particular, when $p=2$, the spatial
profile is continuous but has no positive H\"older exponent at
the origin.
\end{example}

\begin{remark}\label{rem1.5}
(i) Let $n=1$, fix $1<p<2$, and set $\alpha=2/p-1$.
For any $\beta\in(0,\alpha)$, choose
$\tilde p\in(p,2/(1+\beta))$ and define
\begin{equation}\label{1.4}
 b(t,x)=t^{-\frac1{\tilde p}}\mathrm{sign}(x)|x|^\beta,
 \qquad (t,x)\in(0,T]\times\mathbb R.
\end{equation}
Then $b\in L^p([0,T];\mathcal C_l^\beta(\mathbb R))$,
and SDE~\eqref{1.1} with drift~\eqref{1.4}, initial time
$s=0$ and initial condition $X_0=0$ admits two weak solutions
with different laws
(see \cite[Section~1.3]{GG}). Thus, for each fixed $p\in(1,2)$,
weak uniqueness can fail when the spatial exponent is below
$2/p-1$, however close it is to this critical value.
Theorem~\ref{the1.2} reaches the critical exponent
$\alpha=2/p-1$ under an additional Dini-type condition on
the factor $\rho$ in the spatial modulus $r^\alpha\rho(r)$.
The exponent $2p/(5p-2)$ arises from the estimates used in
our proof. This condition is sufficient for the theorem;
its necessity is not established here.

\smallskip
(ii) At the endpoint $p=2$, our theorem and
\cite[Theorem~1.2]{WHLY} impose the same condition on $\rho$,
namely that $\sqrt{\rho}$ be Dini.
The difference lies in the assumptions on the coefficients:
Wei, Hu, Lv and Yuan \cite{WHLY} allow spatially dependent
diffusion coefficients but require spatially bounded drifts,
whereas our theorem considers constant diffusion coefficients
and allows spatially unbounded drifts of at most linear growth.
In the common case of identity diffusion and spatially bounded
drift, both results yield a unique stochastic flow of
homeomorphisms. Our theorem also covers the critical range
$1<p<2$, which is not covered by their theorem.
\end{remark}

\subsection{Main analytic ingredient}\label{sec1.2}
The proof of Theorem~\ref{the1.2} relies on the It\^o--Tanaka trick, which reduces the analysis of SDE \eqref{1.1}
to that of an associated parabolic equation with drift.
This leads naturally to a maximal regularity problem in function spaces adapted to critical H\"older--Dini moduli.
Establishing such a theory is one of the key analytic ingredients of the paper.

\smallskip
To formulate the relevant analytic problem, let ${\mathcal H}$ be a Banach space of functions on ${\mathbb R}^n$, and consider the abstract Cauchy problem
\begin{equation}\label{1.5}
\left\{
\begin{aligned}
&\partial_t u(t,x)=A u(t,x)+f(t,x), \quad (t,x)\in(0,T]\times{\mathbb R}^n,\\
&u(0,x)=0.
\end{aligned}\right.
\end{equation}
We say that $A$ enjoys maximal $L^p_t$--${\mathcal H}_x$ regularity if the solution satisfies
\begin{equation}\label{1.6}
u\in L^p([0,T];{\mathcal D}(A))\cap W^{1,p}([0,T];{\mathcal H})
\quad\text{and}\quad
\|\partial_t u\|_{L^p_t{\mathcal H}_x}
+\|A u\|_{L^p_t{\mathcal H}_x}
\le C\|f\|_{L^p_t{\mathcal H}_x}.
\end{equation}

In the H\"older scale, Krylov \cite{Krylov02} established maximal regularity for $A=\Delta$ on ${\mathcal C}_b^\alpha({\mathbb R}^n)$ with $\alpha\in(0,1)$, obtaining mixed-norm estimates that play a fundamental role in the study of SDEs and stochastic transport equations \cite{FF11,FF12,FF13,FGP,WDGL}.
However, when $\alpha=0$, estimate \eqref{1.6} fails even in the constant-coefficient case, which indicates that the borderline regime requires finer moduli than classical H\"older continuity.

\smallskip
For the elliptic problem
\[
\Delta u(x)=h(x), \qquad x\in B_1=\{x\in \mathbb{R}^n:\ |x|<1\},
\]
Burch \cite{Burch} obtained sharp regularity estimates showing that second derivatives may still be Dini continuous when the modulus of continuity of the data is Dini.
More precisely, if $\rho_h(r):=\sup_{|x-y|<r}|h(x)-h(y)|$ is a two-fold Dini modulus, then $\nabla^2 u$ is again Dini continuous.
Related parabolic extensions were studied in \cite{Wang,WLW}.

\smallskip
These results suggest studying maximal regularity in spaces combining H\"older and Dini moduli.
A natural question is whether such mixed moduli are preserved under the parabolic evolution.

\smallskip
\textbf{(\(\clubsuit\))} Let $\alpha\in(0,1)$, $\beta\in{\mathbb R}$ and $p\in(1,\infty]$.
Suppose $f\in L^p([0,T];{\mathcal C}({\mathbb R}^n))$ satisfies
\begin{align}\label{1.7}
|f(t,x)-f(t,y)|
\le f_1(t)|x-y|^\alpha |\log|x-y||^\beta,
\qquad |x-y|<\tfrac12,
\end{align}
for some $f_1\in L^p([0,T])$.
Does the solution of \eqref{1.5} with $A=\Delta/2$ satisfy
\begin{align}\label{1.8}
|\nabla^2 u(t,x)-\nabla^2 u(t,y)|
\le f_2(t)|x-y|^\alpha |\log|x-y||^\beta, \qquad |x-y|<\tfrac12,
\end{align}
for some $f_2\in L^p([0,T])$ ?

\smallskip
A positive answer to this question provides one of the key analytic ingredients in our treatment of critical SDEs.
To obtain such an answer, we introduce H\"older--Dini type function spaces and establish maximal regularity estimates for non-divergence parabolic equations by combining refined heat kernel bounds with a vector-valued Calder\'on--Zygmund theorem.
These estimates apply both to operators with time-dependent diffusion and to drifted parabolic equations with space-dependent coefficients, thereby providing the analytic framework for the proof of Theorem~\ref{the1.2}.

\smallskip
The rest of the paper is organized as follows.
In Section~\ref{sec2}, we introduce the H\"older--Dini type function spaces and collect preliminary results.
Section~\ref{sec3} is devoted to maximal regularity estimates for parabolic equations in Lebesgue--H\"older--Dini spaces.
Sections~\ref{sec4}--\ref{sec7} contain the proof of Theorem~\ref{the1.2} together with the corresponding maximal regularity results.

\medskip\noindent
\textbf{Notations.} ${\mathbb R}_+=\{r\in{\mathbb R}:\ r>0\}$.
The letter $C$ denotes a positive constant, whose values may change in different places.
${\mathbb N}$ is the set of natural numbers, ${\mathbb N}_0={\mathbb N}\cup \{0\}$ and ${\mathbb Z}$ is the set of integers.

\section{Preliminaries}\label{sec2}\setcounter{equation}{0}
In this section, we introduce the analytic framework used throughout the paper.
These preliminaries are tailored to the maximal regularity estimates needed later
for the proof of Theorem~\ref{the1.2}.
In particular, we collect some Calder\'on--Zygmund type tools and introduce the
H\"older--Dini function spaces adapted to the critical regularity regime of SDE \eqref{1.1}.

\smallskip
We begin with several basic notions.
\begin{definition}\label{def2.1}
Let $\{{\mathbb Q}_k\}_{k\in{\mathbb Z}}$ be a sequence of partitions of $\mathbb{R}$,
where each ${\mathbb Q}_k$ consists of disjoint Borel subsets $Q\subset\mathbb{R}$ and satisfies
\[
R_k:=\sup_{Q\in{\mathbb Q}_k}\operatorname{diam}(Q)<\infty.
\]
We call $\{{\mathbb Q}_k\}_{k\in{\mathbb Z}}$ a filtration of partitions if the following hold:
\begin{enumerate}
\item[(i)] The partitions become finer as $k$ increases; that is,
\[
\inf_{Q\in{\mathbb Q}_k}|Q|\to\infty \ \text{as}\ k\to-\infty,
\qquad
R_k\to0 \ \text{as}\ k\to\infty.
\]
\item[(ii)] The partitions are nested: for each $k$ and each $Q\in{\mathbb Q}_k$,
there exists a unique $Q'\in{\mathbb Q}_{k-1}$ such that $Q\subset Q'$.
\item[(iii)] The regularity property holds:
for $Q$ and $Q'$ as in (ii), we have
\[
|Q'|\le N_0|Q|,
\]
where $N_0>0$ is a constant independent of $k$, $Q$ and $Q'$, and $|\cdot|$ denotes the Lebesgue measure.
\end{enumerate}
\end{definition}

Let ${\mathcal H}$ and $\widetilde{{\mathcal H}}$ be Banach spaces.
We write $L({\mathcal H};\widetilde{{\mathcal H}})$ for the space of bounded linear operators from
${\mathcal H}$ to $\widetilde{{\mathcal H}}$.
Let $D$ be an open or closed domain in $\mathbb{R}$, and let
${\mathcal C}_0^\infty(D;{\mathcal H})$ denote the space of compactly supported smooth
${\mathcal H}$-valued functions on $D$.
We next recall the notion of operator-valued Calder\'on--Zygmund kernels.

\begin{definition}\label{def2.2}
Let $\{{\mathbb Q}_k\}_{k\in{\mathbb Z}}$ be a filtration of partitions.
For each $t,r\in\mathbb{R}$ with $t\neq r$, let
${\mathcal K}(t,r)\in L({\mathcal H};\widetilde{{\mathcal H}})$.
We say that ${\mathcal K}$ is an $L({\mathcal H};\widetilde{{\mathcal H}})$-valued
Calder\'on--Zygmund kernel relative to $\{{\mathbb Q}_k\}_{k\in{\mathbb Z}}$
if the following conditions hold:
\begin{enumerate}
\item[(i)] For every $t\in\mathbb{R}$ and $r_0>0$,
\[
{\mathcal K}(t,\cdot)\in L^1_{\mathrm{loc}}(B_{r_0}^c(t);L({\mathcal H};\widetilde{\mathcal H})),
\]
where $B_{r_0}^c(t):=\{r\in\mathbb{R}:|r-t|\geq r_0\}$.
\item[(ii)] The function
\[
(t,r,\tau)\mapsto
\|{\mathcal K}(t,r)-{\mathcal K}(t,\tau)\|_{L({\mathcal H};\widetilde{{\mathcal H}})}
\]
is measurable for all $(t,r,\tau)\in\mathbb{R}^3$ with $t\neq r$ and $t\neq\tau$.
\item[(iii)] There exists a constant $C_0\geq1$ such that, for each $Q\in{\mathbb Q}_k$, there exists a Borel set $Q^*$ with $\overline{Q}\subset Q^*$, $|Q^*|\le C_0|Q|$, and for every $r\in Q$,
\begin{align}\label{2.1}
\int_{\mathbb{R}\setminus Q^*}
\|{\mathcal K}(t,r)-{\mathcal K}_{|k}(t,r)\|_{L({\mathcal H};\widetilde{{\mathcal H}})}dt
\le C_0,
\end{align}
where
\[
{\mathcal K}_{|k}(t,r)
:=\frac{1}{|Q_k(r)|}\int_{Q_k(r)}{\mathcal K}(t,\tau)d\tau,
\]
and $Q_k(r)$ denotes the unique $\widetilde Q\in{\mathbb Q}_k$ such that $r\in\widetilde Q$.
\end{enumerate}
\end{definition}

\begin{lemma}\label{lem2.3}(\cite[Lemma 1]{TW})
Let ${\mathcal K}$ satisfy conditions (i) and (ii) of Definition~\ref{def2.2}, and let $\{{\mathbb Q}_k\}_{k\in{\mathbb Z}}$ be the filtration of dyadic cubes. Assume that ${\mathcal K}(t,r)$ is weakly differentiable in $r$ for $r\neq t$ and
\[
\|\partial_r{\mathcal K}(t,r)\|_{L({\mathcal H};\widetilde{{\mathcal H}})}\le C\phi(|t-r|)
\]
for all $r\neq t$, with a constant $C>0$ independent of $t,r$ and a function $\phi$ satisfying
\begin{align}\label{2.2}
\iota\int_\iota^{\infty}\phi(\tau)d\tau\le C<\infty
\end{align}
for all $\iota>0$. Then ${\mathcal K}$ is an $L({\mathcal H};\widetilde{{\mathcal H}})$-valued Calder\'on--Zygmund kernel relative to the filtration of dyadic cubes, with the constant $C_0$ in \eqref{2.1} depending only on $C$.
\end{lemma}

The following lemma, due to Krylov, will be used to derive weak- and strong-type estimates from the above Calder\'on--Zygmund framework.

\begin{lemma}\label{lem2.4}(\cite[Theorems 2.5 and 2.9]{Krylov02})
Given a (possibly nonlinear) operator
\[
{\mathcal A}:L^\infty(\mathbb{R};{\mathcal H})\to L^\infty(\mathbb{R}),
\]
suppose that:
\begin{enumerate}
\item[(i)] ${\mathcal A}$ is subadditive and bounded; that is, for some constant $C>0$, every $k\in\mathbb{N}$, and $f,f_m\in L^\infty(\mathbb{R};{\mathcal H})$ $(m=1,\ldots,k)$, we have
\[
\Big|{\mathcal A}\Big(\sum_{m=1}^k f_m(t)\Big)\Big|
\le \sum_{m=1}^k |{\mathcal A}f_m(t)|,\qquad \text{a.e.},
\]
and
\begin{align}\label{2.3}
\|{\mathcal A}f\|_{L^\infty(\mathbb{R})}
\le C\|f\|_{L^\infty(\mathbb{R};{\mathcal H})}.
\end{align}
\item[(ii)] For each $g\in{\mathcal C}_0^\infty(\mathbb{R};{\mathcal H})$ and for almost all $t$ outside the closed support of $g$,
\[
|{\mathcal A}g(t)|
\le \bigg\|\int_{\mathbb{R}}{\mathcal K}(t,r)g(r)dr\bigg\|_{\widetilde{{\mathcal H}}},
\]
where ${\mathcal K}(t,r)$ is an $L({\mathcal H};\widetilde{{\mathcal H}})$-valued Calder\'on--Zygmund kernel relative to a filtration of partitions.
\item[(iii)] If $f,f_1,f_2,\ldots\in L^\infty(\mathbb{R};{\mathcal H})$, $f$ and all $f_k$ vanish outside the same ball, the norms $\|f_k\|_{L^\infty(\mathbb{R};{\mathcal H})}$ are bounded in $k$, and
\[
\|f(t)-f_k(t)\|_{\mathcal H}\to 0
\]
for almost every $t\in\mathbb{R}$, then there exists a subsequence $k(i)\to\infty$ such that
\begin{align}\label{2.4}
|{\mathcal A}f(t)|
\le \liminf_{i\to\infty}|{\mathcal A}f_{k(i)}(t)|,\qquad \text{a.e.}
\end{align}
\end{enumerate}
Then the operator ${\mathcal A}$ is of weak type $(1,1)$ on smooth compactly supported functions; namely, there exists a constant $C_1>0$ such that for any $f\in{\mathcal C}_0^\infty(\mathbb{R};{\mathcal H})$ and $\gamma>0$,
\begin{align}\label{2.5}
\gamma\big|\{t:|{\mathcal A}f(t)|>\gamma\}\big|
\le C_1\int_{\mathbb{R}}|f(t)|_{\mathcal H}dt,
\end{align}
where $|\cdot|$ on the left-hand side denotes the Lebesgue measure of the set $\{t:|{\mathcal A}f(t)|>\gamma\}$. Furthermore, ${\mathcal A}$ is of strong type $(p,p)$ for every $p\in(1,\infty)$; that is, there exists a constant $C_2>0$ such that for all $f\in{\mathcal C}_0^\infty(\mathbb{R};{\mathcal H})$,
\begin{align}\label{2.6}
\|{\mathcal A}f\|_{L^p(\mathbb{R})}
\le C_2\|f\|_{L^p(\mathbb{R};{\mathcal H})}.
\end{align}
\end{lemma}

We next introduce the moduli of continuity that will be used to describe the critical function spaces for the drift and for the associated parabolic equations.

\begin{definition}\label{def2.5}
A continuous, nondecreasing function $\rho:\mathbb{R}_+\to\mathbb{R}_+$ is said to be a Dini function if
\[
\int_0^1\frac{\rho(r)}{r}dr<\infty.
\]
If $\rho$ is a Dini function and additionally satisfies
\[
\int_0^1\int_0^{t_k}\cdots\int_0^{t_2}\rho(t_1)\frac{dt_1}{t_1}\cdots \frac{dt_{k-1}}{t_{k-1}}\frac{dt_k}{t_k}<\infty,
\qquad 0\le t_1\le t_2\le\cdots\le t_k\le 1,
\]
for some integer $k\geq2$, then $\rho$ is said to be a $k$-fold Dini function.
\end{definition}

\begin{definition}\label{def2.6}
A continuous function $\rho:\mathbb{R}_+\to\mathbb{R}_+$ is called slowly varying at zero (in Karamata's sense \cite[p.~6]{BGT}) if for every $\upsilon>0$,
\[
\lim_{r\to0}\frac{\rho(\upsilon r)}{\rho(r)}=1.
\]
\end{definition}

\begin{example}\label{exa2.7}
Let $\beta\in (-\infty,0)$ and define
\[
\rho(r)=
\begin{cases}
|\log r|^\beta, & 0<r<\frac12,\\[0.1cm]
\psi(r), & r\geq \frac12,
\end{cases}
\]
where $\psi\in{\mathcal C}^1(\mathbb{R}_+)$ is chosen such that $\psi'(r)\geq0$ for $r\geq1/2$ and
\[
\psi\Big(\frac12\Big)=|\log2|^\beta,\qquad
\psi'\Big(\frac12\Big)=-2\beta|\log2|^{\beta-1}.
\]
Then $\rho$ is slowly varying at zero but is not a Dini function if $\beta\in [-1,0)$.

\smallskip
If $\beta<-1$, then
\[
\int_0^1\frac{\rho(r)}{r}dr
=\int_0^{1/2}\frac{|\log r|^\beta}{r}dr+\int_{1/2}^1\frac{\psi(r)}{r}dr
\le -\frac{|\log2|^{\beta+1}}{\beta+1}+(\log 2)\psi(1)<\infty.
\]
Moreover, for all $\upsilon>0$,
\[
\lim_{r\to0}\frac{\rho(\upsilon r)}{\rho(r)}
=\lim_{r\to0}\frac{|\log(\upsilon r)|^\beta}{|\log r|^\beta}
=\lim_{r\to0}\frac{|\log\upsilon+\log r|^\beta}{|\log r|^\beta}=1.
\]
Therefore, $\rho$ is both a Dini function and slowly varying at zero.

\smallskip
Similarly, for $\beta\in (-\infty,0)$, define
\[
\rho(r)=
\begin{cases}
|\log r|^\beta e^{-1/r}, & 0<r<\frac12,\\[0.1cm]
\psi_1(r), & r\geq \frac12,
\end{cases}
\]
where $\psi_1\in{\mathcal C}^1(\mathbb{R}_+)$ satisfies $\psi_1'(r)\geq0$ for $r\geq1/2$ and
\[
\psi_1\Big(\frac12\Big)=|\log2|^\beta e^{-2},\qquad
\psi_1'\Big(\frac12\Big)=|\log2|^{\beta-1}e^{-2}\big[-2\beta+4\log2\big].
\]
Then $\rho$ is a Dini function. However, for all $\upsilon>1$,
\[
\lim_{r\to0}\frac{\rho(\upsilon r)}{\rho(r)}
=\lim_{r\to0}\frac{|\log(\upsilon r)|^\beta e^{-1/(\upsilon r)}}{|\log r|^\beta e^{-1/r}}
=\lim_{r\to0}\Big[\frac{|\log\upsilon+\log r|^\beta}{|\log r|^\beta}e^{\frac{\upsilon-1}{\upsilon r}}\Big]
=\infty,
\]
and thus $\rho$ is not slowly varying at zero.
\end{example}

The next lemma gives a useful structural representation of slowly varying moduli, which will be used repeatedly in the regularity analysis below.
\begin{lemma}\label{lem2.8}
Let $\rho: \mathbb R_+\to\mathbb R_+$ be continuous, monotone,
and slowly varying at zero. For each $r_0\in(0,1]$, there
exist continuous functions $c,\zeta:(0,r_0]\to\mathbb R$
such that
\begin{align}\label{2.7}
\rho(r)=\exp\Big\{
c(r)-\int_r^{r_0}\frac{\zeta(\tau)}{\tau}d\tau
\Big\},\qquad 0<r\le r_0,
\end{align}
with
\begin{align}\label{2.8}
\lim_{r\downarrow0}c(r)=c_0\in\mathbb R,
\qquad
\lim_{r\downarrow0}\zeta(r)=0.
\end{align}
Moreover, $\zeta$ can be chosen nonnegative if $\rho$ is
nondecreasing, and nonpositive if $\rho$ is nonincreasing.
In particular,
\begin{align}\label{2.9}
\lim_{r\downarrow0}
\int_r^{r_0}\frac{\zeta(\tau)}{\tau}d\tau
=
\begin{cases}
\ \infty, & \text{if }\rho(r)\to0,\\
-\infty, & \text{if }\rho(r)\to\infty.
\end{cases}
\end{align}
If $\rho$ is constant on $(0,r_0]$, one may take
$c(r)\equiv\log\rho(r_0)$ and $\zeta(r)\equiv0$.
\end{lemma}

\smallskip\noindent
{\bf Proof.}
Fix $r_0\in(0,1]$. For $0<r\le r_0$, define
\[
\zeta(r):=\log\frac{\rho(r)}{\rho(r/e)},
\qquad
c(r):=\log\rho(r)
      +\int_r^{r_0}\frac{\zeta(\tau)}{\tau}d\tau.
\]
The positivity and continuity of $\rho$ imply that $c$ and
$\zeta$ are continuous on $(0,r_0]$. Moreover,
\[
\rho(r)
=\exp\Big\{
c(r)-\int_r^{r_0}\frac{\zeta(\tau)}{\tau}d\tau
\Big\},
\]
which gives \eqref{2.7}. By slow variation of $\rho$ at zero,
\[
\lim_{r\downarrow0}\zeta(r)
=
\lim_{r\downarrow0}
\log\frac{\rho(r)}{\rho(r/e)}
=0.
\]
Thus the second conclusion in \eqref{2.8} holds.

\smallskip
To determine the limit of $c$, a change of variables yields
\[
\begin{aligned}
\int_r^{r_0}\frac{\zeta(\tau)}{\tau}d\tau=
\int_r^{r_0}\frac{\log\rho(\tau)}{\tau}d\tau
-\int_{r/e}^{r_0/e}\frac{\log\rho(\tau)}{\tau}d\tau=
\int_{r_0/e}^{r_0}\frac{\log\rho(\tau)}{\tau}d\tau
-\int_{r/e}^{r}\frac{\log\rho(\tau)}{\tau}d\tau.
\end{aligned}
\]
Consequently,
\[
c(r)
=
\int_{r_0/e}^{r_0}\frac{\log\rho(\tau)}{\tau}d\tau
+\int_{r/e}^{r}
  \frac{\log\rho(r)-\log\rho(\tau)}{\tau}d\tau.
\]
Since $\rho$ is monotone, for $\tau\in[r/e,r]$,
\[
|\log\rho(r)-\log\rho(\tau)|
\le
|\log\rho(r)-\log\rho(r/e)|
=|\zeta(r)|.
\]
It follows that
\[
\bigg|
c(r)-\int_{r_0/e}^{r_0}
          \frac{\log\rho(\tau)}{\tau}d\tau
\bigg|
\le|\zeta(r)|
\longrightarrow0, \quad \text{as} \quad r\downarrow0.
\]
Thus
\[
\lim_{r\downarrow0}c(r)
=
\int_{r_0/e}^{r_0}\frac{\log\rho(\tau)}{\tau}d\tau
=:c_0\in\mathbb R,
\]
which proves the first conclusion in \eqref{2.8}.

\smallskip
If $\rho$ is nondecreasing, then $\rho(r)\geq\rho(r/e)$,
so $\zeta(r)\geq0$. If $\rho$ is nonincreasing, the reverse
inequality gives $\zeta(r)\le0$. Finally,
\[
\int_r^{r_0}\frac{\zeta(\tau)}{\tau}d\tau
=c(r)-\log\rho(r).
\]
Since $c(r)\to c_0\in\mathbb R$, the integral tends to
$\infty$ when $\rho(r)\to0$, and to $-\infty$ when
$\rho(r)\to\infty$. If $\rho$ is constant on $(0,r_0]$,
the above construction gives $\zeta\equiv0$ and
$c\equiv\log\rho(r_0)$.
\qed

\smallskip
We are now in a position to introduce the function spaces that form the basic analytic framework of the paper. These spaces interpolate between classical H\"older classes and Dini-type classes, and are adapted to the critical regime of SDE \eqref{1.1}.

\begin{definition}\label{def2.9}
Let $\alpha\in(0,1)$, let $\rho:\mathbb{R}_+\to\mathbb{R}_+$ be a continuous function, and let $h:\mathbb{R}^n\to\mathbb{R}$ be Borel measurable. Suppose there exists a constant $C>0$ such that
\[
|h(x)-h(y)|\le C|x-y|^\alpha\rho(|x-y|),\qquad x,y\in\mathbb{R}^n,\ \ |x-y|\le1.
\]

\begin{enumerate}
\item[(i)] If $\rho$ is nondecreasing with $\rho(r)\to0$ as $r\to0$, and $r^{-\beta}\rho(r)\to\infty$ as $r\to0$ for every $\beta\in(0,1)$, then $h$ is said to be locally strongly H\"older continuous. The space of all such functions is denoted by ${\mathcal C}_{l,s}^{\alpha,\rho}(\mathbb{R}^n)$. In particular, if $\rho$ is additionally a Dini function, then $h$ is said to be locally H\"older--Dini continuous, and we write ${\mathcal C}_{l,d}^{\alpha,\rho}(\mathbb{R}^n)$ for the corresponding space. If $\rho$ is a $k$-fold Dini function for some integer $k\geq2$, then $h$ is said to be locally H\"older--$k$--Dini continuous, and we denote the space by ${\mathcal C}_{l,d_k}^{\alpha,\rho}(\mathbb{R}^n)$.
\item[(ii)] If $\rho$ is constant on $(0,1]$, then $h$ is said to be locally H\"older continuous. The corresponding space is denoted by ${\mathcal C}_{l,c}^{\alpha,\rho}(\mathbb{R}^n)$.
\item[(iii)] If $\rho$ is nonincreasing with $\rho(r)\to\infty$ as $r\to0$, and $r^\beta\rho(r)\to0$ as $r\to0$ for every $\beta\in(0,1)$, then $h$ is said to be locally weakly H\"older continuous. The corresponding space is denoted by ${\mathcal C}_{l,w}^{\alpha,\rho}(\mathbb{R}^n)$.
\end{enumerate}
\end{definition}

\begin{remark}\label{rem2.10}
(i) From Definition~\ref{def2.9}, the following inclusions hold:
\[
{\mathcal C}_{l,d_m}^{\alpha,\rho}(\mathbb{R}^n)
\subseteq {\mathcal C}_{l,d_k}^{\alpha,\rho}(\mathbb{R}^n)
\subseteq {\mathcal C}_{l,s}^{\alpha,\rho}(\mathbb{R}^n)
\subseteq {\mathcal C}_{l,c}^{\alpha,\rho}(\mathbb{R}^n)
\subseteq {\mathcal C}_{l,w}^{\alpha,\rho}(\mathbb{R}^n),
\qquad 1\le k<m,
\]
where the notation $d_1$ is understood as $d$.
These inclusions reflect the hierarchy of regularity moduli:
Dini conditions impose integrable decay, whereas the strong, classical, and weak H\"older classes successively relax this requirement.

\smallskip
(ii) If $\rho$ is a Dini function, we also use ${\mathcal C}_{l,d}^{0,\rho}(\mathbb{R}^n)$ to denote the set of all continuous functions $h$ on $\mathbb{R}^n$ such that
\[
|h(x)-h(y)|\le C\rho(|x-y|),\qquad x,y\in\mathbb{R}^n,\ \ |x-y|\le1.
\]
\end{remark}

\begin{remark}\label{rem2.11}
\noindent
(i) Let $h\in{\mathcal C}_{l,\varsigma}^{\alpha,\rho}(\mathbb{R}^n)$ with $\varsigma\in\{d_k,s,c,w\}$ ($k\in\mathbb{N}$) and $\alpha\in[0,1)$. For $x,y\in\mathbb{R}^n$ with $|x-y|>1$, there exist points $x_1,x_2,\ldots,x_\eta$ (where $\eta$ is the integer part of $|x-y|$) such that
\[
|x-x_1|=|x_1-x_2|=\cdots=|x_{\eta-1}-x_\eta|=1,\qquad |x_\eta-y|<1.
\]
Set $x_0:=x$. Then
\[|h(x)-h(y)|
\le |h(x_\eta)-h(y)|+\sum_{i=1}^\eta |h(x_{i-1})-h(x_i)|
\le C\big[\rho(1)+\sup_{0<r\le1}r^\alpha\rho(r)\big]|x-y|,
\]
which shows that $h$ grows at most linearly. We define the norm on ${\mathcal C}_{l,\varsigma}^{\alpha,\rho}(\mathbb{R}^n)$ by
\begin{align*}
\|h\|_{{\mathcal C}_{l,\varsigma}^{\alpha,\rho}(\mathbb{R}^n)}
&=\sup_{x\in\mathbb{R}^n}\frac{|h(x)|}{1+|x|}
+\sup_{0<|x-y|\le1}\frac{|h(x)-h(y)|}{|x-y|^\alpha\rho(|x-y|)}\\
&=:\|(1+|\cdot|)^{-1}h(\cdot)\|_0+[h]_{\alpha,\rho}
=: \|h\|_{l,\alpha,\rho}.
\end{align*}
Then ${\mathcal C}_{l,\varsigma}^{\alpha,\rho}(\mathbb{R}^n)$ is a Banach space for every $\varsigma\in\{d_k,s,c,w\}$.

\smallskip
(ii) For a Borel measurable function $h$, if $h$ is bounded and belongs to ${\mathcal C}_{l,\varsigma}^{\alpha,\rho}(\mathbb{R}^n)$, we say that $h\in{\mathcal C}_{b,\varsigma}^{\alpha,\rho}(\mathbb{R}^n)$. For $h\in{\mathcal C}_{b,\varsigma}^{\alpha,\rho}(\mathbb{R}^n)$, we define
\begin{align*}
\|h\|_{{\mathcal C}_{b,\varsigma}^{\alpha,\rho}(\mathbb{R}^n)}
=\sup_{x\in\mathbb{R}^n}|h(x)|
+\sup_{0<|x-y|\le1}\frac{|h(x)-h(y)|}{|x-y|^\alpha\rho(|x-y|)}
=:\|h\|_0+[h]_{\alpha,\rho}
=: \|h\|_{b,\alpha,\rho}.
\end{align*}
Then ${\mathcal C}_{b,\varsigma}^{\alpha,\rho}(\mathbb{R}^n)$ is also a Banach space for every $\varsigma\in\{d_k,s,c,w\}$.

\smallskip
(iii) Let $\rho$ be as in ${\mathcal C}_{\vartheta,\varsigma}^{\alpha,\rho}(\mathbb{R}^n)$ with $\vartheta\in\{l,b\}$, $\varsigma\in\{d_k,s,c,w\}$ and $\alpha\in(0,1)$. Then there exist constants $C_1,C_2>0$ such that for every $r\in(0,1]$ and every $0<\epsilon<\alpha<\beta\le1$,
\begin{align}\label{2.10}
r^\beta\le C_1 r^\alpha\rho(r)\le C_2 r^\epsilon.
\end{align}
\end{remark}

We next pass from pointwise H\"older--Dini regularity to the mixed time-space framework needed for parabolic maximal regularity estimates.
\begin{definition}\label{def2.12}
Let $p\in[1,\infty]$ and $\alpha\in[0,1)$.
Let $\rho:\mathbb{R}_+\to\mathbb{R}_+$ be a continuous modulus and
$\varsigma\in\{s,c,w\}$. We define the locally bounded
Lebesgue--H\"older space
$L^p(\mathbb{R};{\mathcal C}_{l,\varsigma}^{\alpha,\rho}
(\mathbb{R}^n))$ as the set of all Borel measurable functions
$h:\mathbb{R}\times\mathbb{R}^n\to\mathbb{R}$ such that
$\|(1+|\cdot|)^{-1}h(t,\cdot)\|_0\in L^p(\mathbb R)$
and
\[
|h(t,x)-h(t,y)|
\le f(t)|x-y|^\alpha\rho(|x-y|),
\qquad x,y\in\mathbb{R}^n,\quad 0<|x-y|\le1,
\]
for some $f\in L^p(\mathbb{R})$.

\smallskip
If $\rho$ is a $k$-fold Dini function for some $k\geq1$
(with $d_1:=d$), we write
$L^p(\mathbb{R};{\mathcal C}_{l,d_k}^{\alpha,\rho}
(\mathbb{R}^n))$
for the corresponding Lebesgue--H\"older--$k$--Dini space.

\smallskip
Let $m,k\in\mathbb{N}$ and
$\varsigma\in\{d_k,s,c,w\}$.
We say that
$h\in L^p(\mathbb{R};
{\mathcal C}_{l,\varsigma}^{m+\alpha,\rho}(\mathbb{R}^n))$
if $h(t,\cdot)\in{\mathcal C}(\mathbb R^n)$ for almost every $t$,
$\|(1+|\cdot|)^{-1}h(t,\cdot)\|_0\in L^p(\mathbb R)$,
and the following hold:
\begin{itemize}
\item[(i)]
For each $1\le j\le m$ and
$1\le i_1,\ldots,i_j\le n$, $\partial^j_{x_{i_1}\cdots x_{i_j}}h
 \in L^p(\mathbb{R};{\mathcal C}_b(\mathbb{R}^n))$.

\item[(ii)]
For every $1\le i_1,\ldots,i_m\le n$,
$ [\partial^m_{x_{i_1}\cdots x_{i_m}}
 h(t,\cdot)]_{\alpha,\rho}\in L^p(\mathbb{R})$.
\end{itemize}
For
$h\in L^p(\mathbb{R};
{\mathcal C}_{l,\varsigma}^{m+\alpha,\rho}(\mathbb{R}^n))$
with $m\in\mathbb{N}_0$, we define
\begin{align}\label{2.11}
\begin{aligned}
\|h\|_{L^p(\mathbb{R};
{\mathcal C}_{l,\varsigma}^{m+\alpha,\rho}(\mathbb{R}^n))}
&=
\bigg[\int_{\mathbb{R}}
\Big(
\|(1+|\cdot|)^{-1}h(t,\cdot)\|_0
+\sum_{j=1}^m\|\nabla^jh(t,\cdot)\|_0
+[\nabla^mh(t,\cdot)]_{\alpha,\rho}
\Big)^pdt\bigg]^{\frac1p}\\
&=:
\bigg[\int_{\mathbb{R}}
\Big(
\|h(t,\cdot)\|_{l,m,0}
+[\nabla^mh(t,\cdot)]_{\alpha,\rho}
\Big)^pdt\bigg]^{\frac1p}=:
\bigg[\int_{\mathbb{R}}
\|h(t,\cdot)\|_{l,m+\alpha,\rho}^pdt
\bigg]^{\frac1p}.
\end{aligned}
\end{align}
When $p=\infty$, the integral is understood as the
essential supremum.
With this norm, the space
$L^p(\mathbb{R};
{\mathcal C}_{l,\varsigma}^{m+\alpha,\rho}(\mathbb{R}^n))$
is a Banach space.

\smallskip
If $h(t,\cdot)$ is bounded for almost every $t$ and
$\|h(t,\cdot)\|_0\in L^p(\mathbb R)$,
we write
$h\in L^p(\mathbb{R};
{\mathcal C}_{b,\varsigma}^{m+\alpha,\rho}(\mathbb{R}^n))$
and define
\begin{align*}
\|h\|_{L^p(\mathbb{R};
{\mathcal C}_{b,\varsigma}^{m+\alpha,\rho}(\mathbb{R}^n))}
&=
\bigg[\int_{\mathbb{R}}
\Big(
\sum_{j=0}^m\|\nabla^jh(t,\cdot)\|_0
+[\nabla^mh(t,\cdot)]_{\alpha,\rho}
\Big)^pdt\bigg]^{\frac1p}\\
&=:
\bigg[\int_{\mathbb{R}}
\Big(
\|h(t,\cdot)\|_{b,m,0}
+[\nabla^mh(t,\cdot)]_{\alpha,\rho}
\Big)^pdt\bigg]^{\frac1p}=:
\bigg[\int_{\mathbb{R}}
\|h(t,\cdot)\|_{b,m+\alpha,\rho}^pdt
\bigg]^{\frac1p}.
\end{align*}
Finally, for any open or closed interval $D\subset\mathbb{R}$,
we define
$L^p(D;
{\mathcal C}_{\vartheta,\varsigma}^{m+\alpha,\rho}
(\mathbb{R}^n))$
and
$W^{1,p}(D;
{\mathcal C}_{\vartheta,\varsigma}^{m+\alpha,\rho}
(\mathbb{R}^n))$,
with $\vartheta\in\{l,b\}$ and
$\varsigma\in\{d_k,s,c,w\}$, analogously.
\end{definition}

\section{Maximal regularity estimates}\label{sec3}\setcounter{equation}{0}
In this section, we establish the maximal regularity estimates for parabolic equations
in Lebesgue--H\"older--Dini spaces that will be used later in the proof of
Theorem~\ref{the1.2}. We begin with the case of space-independent diffusion coefficients,
and then turn to variable-coefficient equations with drift.

\smallskip
Let $a(t)=(a_{i,j}(t))_{n\times n}$ be a symmetric $n\times n$ matrix-valued, Borel measurable, and bounded function for $t\in \mathbb{R}$.
Assume that there exists a constant $\Gamma\geqq 1$ such that
\begin{align}\label{3.1}
\Gamma^{-1}|\xi|^2\le \sum_{i,j=1}^na_{i,j}(t)\xi_i\xi_j\le \Gamma|\xi|^2,
\end{align}
for all $t\in \mathbb{R}$ and $\xi\in {\mathbb R}^n$.
For $r<t$, set
\[
A_{r,t}:=\int_r^ta(\tau)d\tau,
\qquad
B_{r,t}:=A_{r,t}^{-1}.
\]
Then
\[
\Gamma^{-1}(t-r)|\xi|^2\le \xi^\top A_{r,t}\xi \le \Gamma(t-r)|\xi|^2
\]
and
\[
\Gamma^{-1}(t-r)^{-1}|\xi|^2\le \xi^\top B_{r,t}\xi \le \Gamma(t-r)^{-1}|\xi|^2.
\]
Define
\begin{align}\label{3.2}
K(r,t,x)=1_{\{t>r\}}(2\pi)^{-\frac{n}{2}}\det(B_{r,t})^{\frac{1}{2}}
\exp\Big\{-\tfrac{1}{2}(B_{r,t}x,x)\Big\}
\end{align}
and
\begin{align}\label{3.3}
G_\lambda f(t,x)=\int_{-\infty}^t\int_{{\mathbb R}^n}K(r,t,x-y)f(r,y)e^{-\lambda(t-r)}dydr,
\end{align}
where $\lambda>0$ is fixed and $f\in L^p(\mathbb{R};{\mathcal C}_{\vartheta,\varsigma}^{\alpha,\rho}({\mathbb R}^n))$ with $\vartheta\in \{l,b\}$ and $\varsigma\in \{s,c,w\}$.

\begin{theorem}\label{the3.1}
Let $p\in [1,\infty]$, $\alpha\in (0,1)$ and $\lambda>0$.
Suppose $a(t)=(a_{i,j}(t))_{n\times n}$ is a symmetric, Borel measurable, bounded matrix-valued function such that \eqref{3.1} holds.
Let $f\in L^p(\mathbb{R};{\mathcal C}_{\vartheta,\varsigma}^{\alpha,\rho}({\mathbb R}^n))$ with $\vartheta\in \{l,b\}$ and $\varsigma\in \{s,c,w\}$.
Assume further that $\rho$ is slowly varying at zero. Then there exists a constant $C>0$ such that for every $\gamma>0$:
\begin{enumerate}
\item[(i)] If $p=1$, then
\begin{align}\label{3.4}
\int_{\mathbb{R}}\|G_\lambda f(t,\cdot)\|_{\vartheta,2,0}dt
+\gamma \big|\{t:[\nabla^2 G_\lambda f(t,\cdot)]_{\alpha,\rho}>\gamma\}\big|
\le C\int_{\mathbb{R}}\|f(t,\cdot)\|_{\vartheta,\alpha,\rho}dt.
\end{align}
\item[(ii)] If $p\in (1,\infty]$, then
\begin{align}\label{3.5}
\bigg(\int_{\mathbb{R}}\|G_\lambda f(t,\cdot)\|_{\vartheta,2+\alpha,\rho}^pdt\bigg)^{\frac{1}{p}}
\le C\bigg(\int_{\mathbb{R}}\|f(t,\cdot)\|_{\vartheta,\alpha,\rho}^pdt\bigg)^{\frac{1}{p}},
\end{align}
where the integrals in \eqref{3.5} are understood as essential suprema when $p=\infty$.
\end{enumerate}
\end{theorem}

\begin{remark}\label{rem3.2}
When $p\in[1,\infty]$ and $\alpha>0$, Krylov \cite{Krylov02} proved \eqref{3.4} and \eqref{3.5} for $\vartheta=b$ and $\varsigma=c$.
Theorem~\ref{the3.1} extends Krylov's result not only from bounded H\"older continuous functions to locally H\"older continuous functions, but also from classical H\"older classes to H\"older--Dini type classes.
In particular, it provides the maximal regularity for the second-order differential operator
\[
\frac12\sum_{i,j=1}^n a_{i,j}(t)\partial^2_{x_i,x_j}-\lambda
\]
in the broader framework introduced in Section~\ref{sec2}.
\end{remark}

\begin{example}\label{exa3.3}
For $\beta\in\mathbb{R}$, define
\[
\rho(r)=
\begin{cases}
|\log r|^\beta, & 0<r<\tfrac12,\\[0.1cm]
\psi_1(r)1_{\{\beta<0\}}+\psi_2(r)1_{\{\beta>0\}}+1_{\{\beta=0\}}, & r\geq \tfrac12,
\end{cases}
\]
where $\psi_1,\psi_2$ are smooth functions on $[1/2,\infty)$ such that
\[
\left\{
\begin{array}{ll}
\psi_1(\tfrac{1}{2})=|\log 2|^\beta,\quad
\psi_1'(\tfrac{1}{2})=-2\beta|\log 2|^{\beta-1},\quad
\psi_1'(r)\geq0, & \text{if }\beta<0, \\[0.1cm]
\psi_2(\tfrac{1}{2})=|\log 2|^\beta,\quad
\psi_2'(\tfrac{1}{2})=-2\beta|\log 2|^{\beta-1},\quad
\psi_2'(r)\le0, & \text{if }\beta>0.
\end{array}
\right.
\]
Then $\rho$ is slowly varying at zero.

\smallskip
Let $f\in L^p(\mathbb{R};{\mathcal C}({\mathbb R}^n))$ with $p\in [1,\infty]$.
Suppose there exists $f_1\in L^p(\mathbb{R})$ such that
\[
|f(t,x)-f(t,y)|\le C f_1(t)|x-y|^\alpha \rho(|x-y|),
\qquad x,y\in{\mathbb R}^n,\ \ |x-y|\le1.
\]
Then
\[
\left\{
\begin{array}{ll}
f\in L^p(\mathbb{R};{\mathcal C}_{l,s}^{\alpha,\rho}({\mathbb R}^n)), & \text{if }\beta<0,\\[0.1cm]
f\in L^p(\mathbb{R};{\mathcal C}_{l,c}^{\alpha,\rho}({\mathbb R}^n)), & \text{if }\beta=0,\\[0.1cm]
f\in L^p(\mathbb{R};{\mathcal C}_{l,w}^{\alpha,\rho}({\mathbb R}^n)), & \text{if }\beta>0.
\end{array}
\right.
\]

\smallskip
Let $G_\lambda$ be defined by \eqref{3.3}. Then
\[
\Big[\partial_t-\frac{1}{2}\sum_{i,j=1}^n a_{i,j}(t)\partial^2_{x_i,x_j}+\lambda\Big]G_\lambda f=f,
\]
see, for instance, \cite{Krylov01} when $f\in{\mathcal C}_0^\infty(\mathbb{R}^{n+1})$.
That is, $G_\lambda f$ satisfies
\[
\partial_t u(t,x)=Au(t,x)+f(t,x),
\qquad
A=\frac12\sum_{i,j=1}^n a_{i,j}(t)\partial^2_{x_i,x_j}-\lambda.
\]
By \eqref{3.5}, we have
\[
G_\lambda f\in L^p(\mathbb{R};{\mathcal C}_{\vartheta,\varsigma}^{2+\alpha,\rho}({\mathbb R}^n))
\cap W^{1,p}(\mathbb{R};{\mathcal C}_{\vartheta,\varsigma}^{\alpha,\rho}({\mathbb R}^n)),
\]
and there exists a constant $C>0$ such that
\[
\bigg(\int_{\mathbb{R}}\|\partial_t G_\lambda f(t,\cdot)\|_{\vartheta,\alpha,\rho}^pdt\bigg)^{\frac{1}{p}}
+\bigg(\int_{\mathbb{R}}\|AG_\lambda f(t,\cdot)\|_{\vartheta,\alpha,\rho}^pdt\bigg)^{\frac{1}{p}}
\le C\bigg(\int_{\mathbb{R}}\|f(t,\cdot)\|_{\vartheta,\alpha,\rho}^pdt\bigg)^{\frac{1}{p}},
\]
for $p\in(1,\infty]$.
Therefore, the operator $A$ has maximal regularity for $p\in(1,\infty]$.
\end{example}

Let $T>0$ be fixed. If $f\equiv0$ on $\mathbb{R}\setminus[0,T]$, then by \eqref{3.3} we have
\begin{align}\label{3.6}
G_\lambda f(t,x)=\int_0^t\int_{{\mathbb R}^n}K(r,t,x-y)f(r,y)e^{-\lambda(t-r)}dydr,
\qquad t\in(0,T],
\end{align}
and $\lim_{t\to0}G_\lambda f(t,x)=0$. Hence $G_\lambda f$ satisfies the Cauchy problem
\begin{align}\label{3.7}
\left\{
\begin{aligned}
&\partial_t u(t,x)=\frac12\sum_{i,j=1}^n a_{i,j}(t)\partial^2_{x_i,x_j}u(t,x)-\lambda u(t,x)+f(t,x),
\quad (t,x)\in(0,T]\times{\mathbb R}^n,\\
&u(0,x)=0,\quad x\in{\mathbb R}^n.
\end{aligned}
\right.
\end{align}
Moreover, by the calculations in Section~\ref{sec4}, inequality \eqref{3.5} remains valid for all $\lambda\geq0$ since $f$ vanishes outside $[0,T]$. In particular, \eqref{1.8} holds whenever \eqref{1.7} is satisfied, thereby giving a positive answer to question $(\clubsuit)$. Furthermore, the strong solution of \eqref{3.7} in the class
\[
L^p([0,T];{\mathcal C}_{\vartheta,\varsigma}^{2+\alpha,\rho}({\mathbb R}^n))
\cap W^{1,p}([0,T];{\mathcal C}_{\vartheta,\varsigma}^{\alpha,\rho}({\mathbb R}^n)),
\qquad \vartheta\in\{l,b\},\ \varsigma\in\{d,s,c,w\},
\]
is unique for $p\in(1,\infty]$. Here the unknown function $u(t,x)$ is called a strong solution of \eqref{3.7} if
\[
u,\ \partial_tu,\ \partial^2_{x_i,x_j}u\in L^1([0,T];L^\infty_{\mathrm{loc}}({\mathbb R}^n)),
\qquad 1\le i,j\le n,
\]
which grow at most linearly, such that \eqref{3.7} holds for almost all $(t,x)\in[0,T]\times{\mathbb R}^n$.

\smallskip
For the unique strong solution of \eqref{3.7}, we also obtain boundedness of its gradient.

\begin{theorem}\label{the3.4}
Let $\alpha,p,\rho,a$ and $\vartheta$ be as in Theorem~\ref{the3.1} with $\alpha\geq 2/p-1$, and let $\varsigma\in\{d,s,c,w\}$. Suppose $T>0$, $f\in L^p([0,T];{\mathcal C}_{\vartheta,\varsigma}^{\alpha,\rho}({\mathbb R}^n))$ and $\lambda\geq0$. Let $u$ be defined by \eqref{3.6}.
\begin{enumerate}
\item[(i)] If $\alpha\notin\{2/p-1,2/p\}$, then $u\in L^\infty([0,T];{\mathcal C}_{\vartheta,\varsigma}^{2+\alpha-\frac{2}{p},\rho}({\mathbb R}^n))$
and there exists a constant $C>0$ such that
\begin{align}\label{3.8}
\|u\|_{L^\infty([0,T];{\mathcal C}_{\vartheta,\varsigma}^{2+\alpha-\frac{2}{p},\rho}({\mathbb R}^n))}
\le C\|f\|_{L^p([0,T];{\mathcal C}_{\vartheta,\varsigma}^{\alpha,\rho}({\mathbb R}^n))}.
\end{align}
\item[(ii)] Let $\varsigma=d$.
If $\alpha=2/p-1$, then $u\in L^\infty([0,T];{\mathcal C}_{\vartheta,s}^{1,\hat{\rho}}({\mathbb R}^n))$; and if $\alpha=2/p$, then $u\in L^\infty([0,T];{\mathcal C}_{\vartheta,s}^{2,\hat{\rho}}({\mathbb R}^n))$.
Moreover, there exists a constant $C>0$ such that
\begin{align}\label{3.9}
\|u\|_{L^\infty([0,T];{\mathcal C}_{\vartheta,s}^{2+\alpha-\frac{2}{p},\hat{\rho}}({\mathbb R}^n))}
\le C\|f\|_{L^p([0,T];{\mathcal C}_{\vartheta,d}^{\alpha,\rho}({\mathbb R}^n))},
\qquad \alpha\in\Big\{\tfrac{2}{p}-1,\tfrac{2}{p}\Big\},
\end{align}
where
\begin{align}\label{3.10}
\hat{\rho}(r)=\int_0^r\frac{\rho(\tau)}{\tau}d\tau+\rho(r)+r\int_r^1\frac{\rho(\tau)}{\tau^2}d\tau,
\qquad r\in(0,1].
\end{align}
\end{enumerate}
\end{theorem}

\begin{remark}\label{rem3.5}
(i) Since \eqref{3.5} holds for all $\lambda\geq0$, by \eqref{3.6} and choosing $\varsigma=d$, we obtain
\[
\|u\|_{L^p([0,T];{\mathcal C}_{\vartheta,d}^{2+\alpha,\rho}({\mathbb R}^n))}
\le C\|f\|_{L^p([0,T];{\mathcal C}_{\vartheta,d}^{\alpha,\rho}({\mathbb R}^n))},
\qquad p\in(1,\infty].
\]
In particular, if $f\in L^p([0,T];{\mathcal C}_{b,d_k}^{\alpha,\rho}({\mathbb R}^n))$ with $k\geq1$, then $u\in L^p([0,T];{\mathcal C}_{b,d_k}^{2+\alpha,\rho}({\mathbb R}^n))$.
Thus maximal Lebesgue--H\"older--$k$--Dini regularity is also preserved.

\smallskip
(ii) From \eqref{3.9} and \eqref{3.10}, let $f\in L^p([0,T];{\mathcal C}_{b,d_k}^{\alpha,\rho}({\mathbb R}^n))$ with $k\geq1$. Then the following hold:

\smallskip
$\bullet$ If $k=1$, the modulus $\hat\rho$ need not be Dini.
Thus the estimate does not, in general, guarantee Dini
continuity of $\nabla u$ when $\alpha=2/p-1$, or of
$\nabla^2u$ when $\alpha=2/p$.

\smallskip
$\bullet$ If $k>1$, then
\[
u\in L^\infty\!\big([0,T];{\mathcal C}_{\vartheta,d_{k-1}}^{2+\alpha-\frac{2}{p},\hat{\rho}}({\mathbb R}^n)\big),
\qquad \alpha\in\Big\{\tfrac{2}{p}-1,\tfrac{2}{p}\Big\}.
\]
This refined estimate is consistent with, and further extends, the classical $L^\infty$-theory for elliptic equations \cite{Burch} and parabolic equations \cite{Wang}.
\end{remark}

We now extend the Cauchy problem \eqref{3.7} from the case of space-independent diffusion without drift to the case of space-dependent diffusion with drift, and establish the unique strong solvability. Specifically, we consider
\begin{align}\label{3.11}
\left\{
\begin{aligned}
&\partial_tu(t,x)=\frac12\sum_{i,j=1}^n a_{i,j}(t,x)\partial^2_{x_i,x_j}u(t,x)
+g(t,x)\cdot\nabla u(t,x)\\
&\qquad\qquad
-\lambda u(t,x)+f(t,x),\qquad (t,x)\in(0,T]\times{\mathbb R}^n,\\
&u(0,x)=0,\qquad x\in{\mathbb R}^n,
\end{aligned}
\right.
\end{align}
where $a_{i,j}(t,x)$, $i,j=1,\ldots,n$, are real-valued functions such that
\[
a_{i,j}\in L^\infty([0,T];{\mathcal C}_{b,\varsigma}^{\alpha,\rho}({\mathbb R}^n)),
\qquad \varsigma\in\{d_k,s,c,w\},\ \ k\in\mathbb{N}.
\]
The notion of strong solution for \eqref{3.11} is the same as that for \eqref{3.7}.

\begin{theorem}\label{the3.6}
(\emph{Maximal Lebesgue--H\"older--Dini and Lebesgue--H\"older regularity})
Let the parameters $\alpha,p,\rho,\vartheta,\lambda$ and $f$ be as in Theorem~\ref{the3.4}. Suppose $k\in\mathbb{N}$ and $\varsigma\in\{d_k,s,c,w\}$. Set
\begin{align}\label{3.12}
{\mathcal X}:=L^p([0,T];{\mathcal C}_{\vartheta,\varsigma}^{2+\alpha,\rho}({\mathbb R}^n))
\cap W^{1,p}([0,T];{\mathcal C}_{\vartheta,\varsigma}^{\alpha,\rho}({\mathbb R}^n)).
\end{align}
Let $g=(g_1,\ldots,g_n)\in L^q([0,T];{\mathcal C}_{b,\varsigma}^{\alpha,\rho}({\mathbb R}^n;{\mathbb R}^n))$.
Assume further that $a_{i,j}\in L^\infty([0,T];{\mathcal C}_{b,\varsigma}^{\alpha,\rho}({\mathbb R}^n))$, $1\le i,j\le n$, and that there exists a constant $\Gamma\geq1$ such that
\[
\Gamma^{-1}|\xi|^2\le\sum_{i,j=1}^n a_{i,j}(t,x)\xi_i\xi_j\le\Gamma|\xi|^2,
\qquad \forall\,\xi\in{\mathbb R}^n,\ (t,x)\in[0,T]\times{\mathbb R}^n.
\]
\begin{enumerate}
\item[(i)] If $q=p\in[2,\infty]$, then there exists a unique solution $u\in{\mathcal X}$ to \eqref{3.11}, satisfying
\begin{align}\label{3.13}
\|\partial_tu\|_{L^p([0,T];{\mathcal C}_{\vartheta,\varsigma}^{\alpha,\rho}({\mathbb R}^n))}
+\|u\|_{L^p([0,T];{\mathcal C}_{\vartheta,\varsigma}^{2+\alpha,\rho}({\mathbb R}^n))}
\le C\|f\|_{L^p([0,T];{\mathcal C}_{\vartheta,\varsigma}^{\alpha,\rho}({\mathbb R}^n))}
\end{align}
for some constant $C>0$. Moreover:

\noindent
$\bullet$ If $\alpha\neq 2/p$, then $u\in L^\infty([0,T];{\mathcal C}_{\vartheta,\varsigma}^{2+\alpha-\tfrac{2}{p},\rho}({\mathbb R}^n))$ and \eqref{3.8} holds.

\noindent
$\bullet$ If $\alpha=2/p$ and $\varsigma=d_k$, then $u\in L^\infty([0,T];{\mathcal C}_{\vartheta,\widetilde\varsigma}^{2,\hat\rho}({\mathbb R}^n))$
and estimate \eqref{3.9} holds, where $\widetilde\varsigma=s$ if $k=1$, and $\widetilde\varsigma=d_{k-1}$ if $k>1$.
\item[(ii)] If $q=\infty$ and $p\in(1,\infty]$, then the same conclusions remain valid: the problem \eqref{3.11} has a unique solution $u\in{\mathcal X}$ satisfying \eqref{3.13}. Furthermore:

\noindent
$\bullet$ If $2/p-1<\alpha\neq 2/p$, then $u\in L^\infty([0,T];{\mathcal C}_{\vartheta,\varsigma}^{2+\alpha-\tfrac{2}{p},\rho}({\mathbb R}^n))$ and \eqref{3.8} holds.

\noindent
$\bullet$ If $\alpha\in\{2/p-1,2/p\}$ and $\varsigma=d_k$, then $u\in L^\infty([0,T];{\mathcal C}_{\vartheta,\widetilde\varsigma}^{2+\alpha-\tfrac{2}{p},\hat\rho}({\mathbb R}^n))$
and \eqref{3.9} holds, where $\widetilde\varsigma=s$ if $k=1$, and $\widetilde\varsigma=d_{k-1}$ if $k>1$.
\end{enumerate}
\end{theorem}

\begin{remark}\label{rem3.7} The maximal regularity results obtained in Theorem \ref{the3.6} provide the analytic input for the It\^o--Tanaka trick used in Section 7. In particular, they will be the key ingredient in the proof of the stochastic flow result stated in Theorem~\ref{the1.2}.
\end{remark}

\section{Proof of Theorem \ref{the3.1}}\label{sec4}
\setcounter{equation}{0}
The proof for $\vartheta=b$ is similar to that for $\vartheta=l$, so we only give the detailed calculation for $\vartheta=l$. Let $G_\lambda f(t,x)$ be given by \eqref{3.3}. Then
\begin{align*}
\sup_{x\in {\mathbb R}^n}\frac{|G_\lambda f(t,x)|}{1+|x|}
&\le C\sup_{x\in {\mathbb R}^n}\frac{1}{1+|x|}\int_{-\infty}^t\int_{{\mathbb R}^n}K(r,t,y)\big[1+|x|+|y|\big]f_1(r)e^{-\lambda(t-r)}dydr \\
&\le C\int_{-\infty}^t f_1(r)e^{-\lambda(t-r)}dr
 +C\int_{-\infty}^t\int_{{\mathbb R}^n}K(r,t,y)|y|f_1(r)e^{-\lambda(t-r)}dydr \\
&=C\int_{-\infty}^t f_1(r)\big[1+(t-r)^{\frac12}\big]e^{-\lambda(t-r)}dr,
\end{align*}
where $f_1(r)=\|(1+|\cdot|)^{-1}f(r,\cdot)\|_0\in L^p(\mathbb{R})$. By Young's inequality,
\begin{align}\label{4.1}
\bigg(\int_{\mathbb{R}}\|(1+|\cdot|)^{-1}G_\lambda f(t,\cdot)\|_0^pdt\bigg)^{\frac{1}{p}}
\le C\bigg(\int_{\mathbb{R}}|f_1(t)|^pdt\bigg)^{\frac{1}{p}},
\end{align}
where the integrals in \eqref{4.1} are understood as the essential supremum when $p=\infty$.

\smallskip
For $1\le i\le n$,
\begin{align}\label{4.2}
\begin{aligned}
|\partial_{x_i}G_\lambda f(t,x)|
&\le 2\int_{-\infty}^t\int_{{\mathbb R}^n}|\partial_{y_i}K(r,t,y)|f_2(r)\big(\mathbf{1}_{\{|y|\le 1\}}|y|^{\alpha}\rho(|y|)+|y|\mathbf{1}_{\{|y|>1\}}\big)e^{-\lambda(t-r)}dydr\\
&\le C\int_{-\infty}^t\int_{{\mathbb R}^n}(t-r)^{-\frac n2}\exp\Big(-\frac{|y|^2}{C\Gamma(t-r)}\Big)f_2(r)\Big[1+|y|^{\frac{\alpha}{2}-1}\Big]e^{-\lambda(t-r)}dydr\\
&\le C\int_{-\infty}^t f_2(r)\Big[1+(t-r)^{\frac{\alpha-2}{4}}\Big]e^{-\lambda(t-r)}dr,
\end{aligned}
\end{align}
where $f_2(r)=[f(r,\cdot)]_{\alpha,\rho}\in L^p(\mathbb{R})$, and in the third line we used \eqref{2.10}. Hence,
\begin{align}\label{4.3}
\bigg(\int_{\mathbb{R}}\|\partial_{x_i} G_\lambda f(t,\cdot)\|_0^pdt\bigg)^{\frac{1}{p}}
\le C\bigg(\int_{\mathbb{R}}[f(t,\cdot)]_{\alpha,\rho}^pdt\bigg)^{\frac{1}{p}}.
\end{align}
Furthermore, for $1\le i,j\le n$, we similarly obtain
\[
|\partial^2_{x_i,x_j}G_\lambda f(t,x)|
\le C \int_{-\infty}^t f_2(r)\Big[(t-r)^{-\frac12}+(t-r)^{\frac{\alpha}{4}-1}\Big]e^{-\lambda(t-r)}dr,
\]
which implies
\begin{align}\label{4.4}
\bigg(\int_{\mathbb{R}}\|\partial^2_{x_i,x_j} G_\lambda f(t,\cdot)\|_0^pdt\bigg)^{\frac{1}{p}}
\le C\bigg(\int_{\mathbb{R}}[f(t,\cdot)]_{\alpha,\rho}^pdt\bigg)^{\frac{1}{p}}.
\end{align}
Combining \eqref{4.1}, \eqref{4.3} and \eqref{4.4}, it remains to show, for every $\gamma>0$,
\begin{align}\label{4.5}
\gamma\big|\{t:\,[\nabla^2 G_\lambda f(t,\cdot)]_{\alpha,\rho}>\gamma\}\big|
\le C\int_{\mathbb{R}}\|f(t,\cdot)\|_{l,\alpha,\rho}dt, \quad {\text if} \ \ p=1,
\end{align}
and
\begin{align}\label{4.6}
\bigg(\int_{\mathbb{R}}[\nabla^2G_\lambda f(t,\cdot)]_{\alpha,\rho}^pdt\bigg)^{\frac{1}{p}}
\le C\bigg(\int_{\mathbb{R}}\|f(t,\cdot)\|_{l,\alpha,\rho}^pdt\bigg)^{\frac{1}{p}}, \quad {\text if} \ \ p\in (1,\infty].
\end{align}

We first prove \eqref{4.6} for $p=\infty$. By \eqref{4.4}, it suffices to show that for every $1\le i,j\le n$ and every $x,y\in {\mathbb R}^n$ with $|x-y|\le 1/3$, there exists $C>0$ such that
\begin{align}\label{4.7}
|\partial^2_{x_i,x_j}G_\lambda f(t,x)-\partial^2_{y_i,y_j}G_\lambda f(t,y)|
\le C\operatorname*{esssup}_{r\in \mathbb{R}}[f(r,\cdot)]_{\alpha,\rho} |x-y|^\alpha\rho(|x-y|), \quad t\in \mathbb{R}.
\end{align}

By \eqref{3.3},
\begin{align}\label{4.8}
&\partial^2_{x_i,x_j}G_\lambda f(t,x)-\partial^2_{y_i,y_j}G_\lambda f(t,y)\notag\\
&=:J_1(t,x,y)+J_2(t,x,y)+J_3(t,x,y)+J_4(t,x,y)+J_5(t,x,y),
\end{align}
where
\begin{align*}
\left\{
\begin{array}{ll}J_1(t,x,y)=\int_{-\infty}^t e^{-\lambda(t-r)}dr\int_{|x-z|\le 2|x-y|}\partial^2_{x_i,x_j}K(r,t,x-z)[f(r,z)-f(r,x)]dz,  \\ [0.2cm]J_2(t,x,y)=-\int_{-\infty}^t e^{-\lambda(t-r)}dr\int_{|x-z|\le 2|x-y|}\partial^2_{y_i,y_j}K(r,t,y-z)[f(r,z)-f(r,y)]dz\\ [0.2cm]
J_3(t,x,y)=\int_{-\infty}^t e^{-\lambda(t-r)}\!dr\int_{|x-z|> 2|x-y|}\partial^2_{y_i,y_j}K(r,t,y-z)[f(r,y)-f(r,x)]dz, \\ [0.2cm] J_4(t,x,y)=\int_{-\infty}^t e^{-\lambda(t-r)}\!dr\int_{2|x-y|<|x-z|\le 1 }\big[\partial^2_{x_i,x_j}K(r,t,x-z)-\partial^2_{y_i,y_j}K(r,t,y-z)\big]
\\ [0.2cm]
\qquad\qquad \quad\times
[f(r,z)-f(r,x)]dz, \\ [0.2cm]J_5(t,x,y)=\int_{-\infty}^t e^{-\lambda(t-r)}dr\int_{|x-z|>1}\big[\partial^2_{x_i,x_j}K(r,t,x-z)
-\partial^2_{y_i,y_j}K(r,t,y-z)\big]\\ [0.2cm]
\qquad\qquad \quad\times[f(r,z)-f(r,x)]dz
\end{array}\right.
\end{align*}

For $J_1$,
\begin{align}\label{4.9}
\begin{aligned}
|J_1(t,x,y)|
&\le C\operatorname*{esssup}_{r\in \mathbb{R}}[f(r,\cdot)]_{\alpha,\rho}
\int_{|z|\le 2|x-y|} |z|^\alpha\rho(|z|)dz \int_0^{\infty} r^{-\frac{n+2}{2}}e^{-\frac{|z|^2}{4\Gamma r}}dr\\
&\le C\operatorname*{esssup}_{r\in \mathbb{R}}[f(r,\cdot)]_{\alpha,\rho}\int_0^{2|x-y|} \frac{\rho(r)}{r^{1-\alpha}}dr.
\end{aligned}
\end{align}
If $|x-z|\le 2|x-y|$, then $|y-z|\le 3|x-y|$, and thus
\[
|J_2(t,x,y)|\le \int_{-\infty}^{t} e^{-\lambda(t-r)}dr \int_{|y-z|\le 3|x-y|}
\big|\partial^2_{y_i,y_j}K(r,t,y-z)\big||f(r,z)-f(r,y)|dz,
\]
which yields
\begin{align}\label{4.10}
|J_2(t,x,y)|\le C\operatorname*{esssup}_{r\in \mathbb{R}}[f(r,\cdot)]_{\alpha,\rho}
\int_0^{3|x-y|} \frac{\rho(r)}{r^{1-\alpha}}dr.
\end{align}

For $J_3$, by the Gauss--Green formula,
\begin{align*}
|J_3(t,x,y)|
&=\bigg|\int_{-\infty}^t e^{-\lambda(t-r)}dr\int_{|x-z|=2|x-y|}
\partial_{y_j}K(r,t,y-z)\nu_i[f(r,y)-f(r,x)]dS\bigg| \\
&\le C\int_{-\infty}^t e^{-\lambda(t-r)}dr \int_{|x-z|=2|x-y|}
|y-z|(t-r)^{-\frac{n+2}{2}}e^{-\frac{|y-z|^2}{2\Gamma(t-r)}}|f(r,y)-f(r,x)|dS,
\end{align*}
where $\nu=(\nu_1,\ldots,\nu_n)$ is the outer unit normal to $\{z:|x-z|=2|x-y|\}$. Hence
\begin{align}\label{4.11}
|J_3(t,x,y)|\le C\operatorname*{esssup}_{r\in \mathbb{R}}[f(r,\cdot)]_{\alpha,\rho}
|x-y|^\alpha\rho(|x-y|).
\end{align}

For $J_4$, since $|x-z|>2|x-y|$, for every $\xi$ on the segment $[x,y]$,
\[
\tfrac12|x-z|\le |\xi-z|\le 2|x-z|.
\]
By the mean value inequality,
\begin{align}\label{4.12}
\begin{aligned}
|J_4(t,x,y)|
&\le C|x-y|\int_{-\infty}^tdr \int_{2|x-y|<|x-z|\le 1}\!|f(r,z)-f(r,x)|
(t-r)^{-\frac{n+3}{2}}e^{-\frac{|x-z|^2}{16\Gamma(t-r)}}dz\\
&\le C|x-y|\operatorname*{esssup}_{r\in \mathbb{R}}[f(r,\cdot)]_{\alpha,\rho}
\int_{2|x-y|<|z|\le 1}|z|^\alpha\rho(|z|)|z|^{-n-1}dz \int_0^{\infty} r^{\frac{n-1}{2}}e^{-\frac{r}{16}}dr\\
&\le C\operatorname*{esssup}_{r\in \mathbb{R}}[f(r,\cdot)]_{\alpha,\rho}
|x-y|\int_{2|x-y|}^1\frac{\rho(r)}{r^{2-\alpha}}dr.
\end{aligned}
\end{align}

Fix $\beta\in(\alpha,1)$. For $|x-z|>2|x-y|$,
\[
\big|\partial^2_{x_i,x_j}K(r,t,x-z)-\partial^2_{y_i,y_j}K(r,t,y-z)\big|
\le C|x-y|^\beta(t-r)^{-\frac{n+2+\beta}{2}}e^{-\frac{|x-z|^2}{16\Gamma(t-r)}},
\]
and hence
\begin{align}\label{4.13}
\begin{aligned}
|J_5(t,x,y)|
&\le C|x-y|^\beta\int_{-\infty}^t e^{-\lambda(t-r)}dr
\int_{|x-z|>1}|f(r,z)-f(r,x)|(t-r)^{-\frac{n+2+\beta}{2}}e^{-\frac{|x-z|^2}{16\Gamma(t-r)}}dz\\
&\le  C|x-y|^\beta\operatorname*{esssup}_{r\in \mathbb{R}}[f(r,\cdot)]_{\alpha,\rho}
\int_0^{\infty}e^{-\lambda r} r^{-\frac{1+\beta}{2}}dr \\
&\le  C|x-y|^\beta\operatorname*{esssup}_{r\in \mathbb{R}}[f(r,\cdot)]_{\alpha,\rho},
\end{aligned}
\end{align}
where, in the second inequality, we used that $f(r,\cdot)$ has at most linear growth (see Remark \ref{rem2.11} (i)).

\smallskip
Combining \eqref{4.8}--\eqref{4.13} gives
\begin{align}\label{4.14}
\begin{aligned}
&|\partial^2_{x_i,x_j}G_\lambda f(t,x)-\partial^2_{y_i,y_j}G_\lambda f(t,y)|\\
&\le C\operatorname*{esssup}_{r\in \mathbb{R}}[f(r,\cdot)]_{\alpha,\rho}
\bigg[\int_0^{3|x-y|} \frac{\rho(r)}{r^{1-\alpha}}dr+|x-y|^\alpha\rho(|x-y|)\\ &\qquad +|x-y|\int_{2|x-y|}^1\frac{\rho(r)}{r^{2-\alpha}}dr+|x-y|^\beta\bigg]\\
&\le C\bigg[\int_0^{|x-y|}\frac{\rho(r)}{r^{1-\alpha}}dr + |x-y|^\alpha\rho(|x-y|)
 +|x-y|\int_{|x-y|}^1\frac{\rho(r)}{r^{2-\alpha}}dr
 +\int_{|x-y|}^{3|x-y|}\frac{\rho(r)}{r^{1-\alpha}}dr\bigg],
\end{aligned}
\end{align}
where in the second inequality we have used \eqref{2.10}.

\smallskip
By L'Hospital's rule, we get
\begin{align*}
\lim_{|x-y|\rightarrow 0}\frac{\int_{|x-y|}^{3|x-y|} \frac{\rho(r)}{r^{1-\alpha}} dr}{\int_0^{|x-y|} \frac{\rho(r)}{r^{1-\alpha}} dr}=\lim_{|x-y|\rightarrow 0}\frac{ \frac{3\rho(3|x-y|)}{(3|x-y|)^{1-\alpha}}-\frac{\rho(|x-y|)}{(|x-y|)^{1-\alpha}}}
{\frac{\rho(|x-y|)}{|x-y|^{1-\alpha}}} = 3^\alpha-1>0,
\end{align*}
where the assumption that $\rho$ is slowly varying at zero has been used. Hence, there is a constant $C>0$ such that for every $|x-y|\le 1/3$
\begin{align}\label{4.15}
\int_{|x-y|}^{3|x-y|} \frac{\rho(r)}{r^{1-\alpha}} dr\le C\int_0^{|x-y|} \frac{\rho(r)}{r^{1-\alpha}} dr.
\end{align}

Substituting \eqref{4.15} into \eqref{4.14} yields
\begin{align}\label{4.16}
&|\partial^2_{x_i,x_j}G_\lambda f(t,x)-\partial^2_{y_i,y_j}G_\lambda f(t,y)|
\notag\\ &\le C\bigg[\int_0^{|x-y|}\frac{\rho(r)}{r^{1-\alpha}}dr+ |x-y|^\alpha\rho(|x-y|)
+|x-y|\int_{|x-y|}^1\frac{\rho(r)}{r^{2-\alpha}}dr\bigg].
\end{align}

If $f\in L^\infty(\mathbb{R};{\mathcal C}_{l,\varsigma}^{\alpha,\rho}({\mathbb R}^n))$ with $\varsigma=s$, then by \eqref{2.7}--\eqref{2.9} and L'Hospital's rule,
\begin{align*}
\lim_{|x-y|\to 0}\frac{\int_0^{|x-y|}\frac{\rho(r)}{r^{1-\alpha}}dr}{|x-y|^\alpha\rho(|x-y|)}
&= e^{-c_0}\lim_{|x-y|\to 0}\frac{\int_0^{|x-y|}\frac{\rho(r)}{r^{1-\alpha}}dr}{|x-y|^\alpha\rho(|x-y|)e^{-c(|x-y|)}} \notag\\
&= e^{-c_0}\lim_{|x-y|\to 0}\frac{\rho(|x-y|)|x-y|^{\alpha-1}}{[\alpha+\zeta(|x-y|)]|x-y|^{\alpha-1}\rho(|x-y|)e^{-c(|x-y|)}}
=\frac{1}{\alpha},
\end{align*}
and
\begin{align*}
\lim_{|x-y|\to 0}\frac{|x-y|\int_{|x-y|}^1\frac{\rho(r)}{r^{2-\alpha}}dr}{|x-y|^\alpha\rho(|x-y|)}
&=e^{-c_0}\lim_{|x-y|\to 0}\frac{-\rho(|x-y|)|x-y|^{\alpha-2}}{[\alpha-1+\zeta(|x-y|)]|x-y|^{\alpha-2}\rho(|x-y|)e^{-c(|x-y|)}}
\notag\\
&=\frac{1}{1-\alpha}.
\end{align*}
Therefore, whenever $|x-y|\le 1/3$,
\begin{align}\label{4.17}
\max\bigg\{\int_0^{|x-y|}\frac{\rho(r)}{r^{1-\alpha}}dr,\ \ |x-y|\int_{|x-y|}^1\frac{\rho(r)}{r^{2-\alpha}}dr\bigg\}
\le C|x-y|^\alpha\rho(|x-y|).
\end{align}

For $\varsigma=c$, inequality \eqref{4.17} holds directly, while for $\varsigma=w$,
\begin{align}\label{4.18}
|x-y|\displaystyle\int_{|x-y|}^1\frac{\rho(r)}{r^{2-\alpha}}dr\le |x-y|\rho(|x-y|)\displaystyle\int_{|x-y|}^1\frac{1}{r^{2-\alpha}}dr\le
\frac{1}{1-\alpha}|x-y|^\alpha\rho(|x-y|),
\end{align}
since $\rho$ is nonincreasing.

\smallskip
Moreover, applying L'Hospital's rule and noting that $\rho$ is slowly varying at zero, it follows that
\begin{align}\label{4.19}
\lim_{|x-y|\to 0}\frac{\int_0^{|x-y|}\frac{\rho(r)}{r^{1-\alpha}}dr}{\int_{|x-y|}^{2|x-y|}\frac{\rho(r)}{r^{1-\alpha}}dr}
&= \lim_{|x-y|\to 0}\frac{\rho(|x-y|)}{2^\alpha\rho(2|x-y|)-\rho(|x-y|)} =\frac{1}{2^\alpha-1}.
\end{align}
Consequently, there exists a constant $C>0$ such that
\begin{align}\label{4.20}
\int_0^{|x-y|}\frac{\rho(r)}{r^{1-\alpha}}dr\le C\int_{|x-y|}^{2|x-y|}\frac{\rho(r)}{r^{1-\alpha}}dr\le C|x-y|^\alpha\rho(|x-y|).
\end{align}
Combining \eqref{4.19} and \eqref{4.20} confirms that \eqref{4.17} remains valid. Hence, \eqref{4.7} follows from \eqref{4.16} and \eqref{4.17}.

\medskip
Let ${\mathcal H}={\mathcal C}_{l,\varsigma}^{\alpha,\rho}({\mathbb R}^n)$ and $\widetilde{{\mathcal H}}={\mathcal C}_{b,\varsigma}^{\alpha,\rho}({\mathbb R}^n)$, $\varsigma\in \{s,c,w\}$. For $1\le i,j\le n$ and $h\in {\mathcal H}$, set
\begin{align}\label{4.21}
{\mathcal K}_{i,j}(t,r)h(x)=\int_{{\mathbb R}^n}\partial^2_{x_i,x_j}K(r,t,x-y)h(y)dy
=:\int_{{\mathbb R}^n}K_{i,j}(r,t,x-y)h(y)dy.
\end{align}
For $f\in L^\infty(\mathbb{R};{\mathcal H})$, define
\begin{align}\label{4.22}
{\mathcal A}_{i,j}^\lambda f(t,\cdot)=\int_{\mathbb{R}}e^{-\lambda(t-r)}{\mathcal K}_{i,j}(t,r)f(r,\cdot)dr
\end{align}
and
\begin{align}\label{4.23}
{\mathcal A} f(t)=\|{\mathcal A}_{i,j}^\lambda f(t,\cdot)\|_{\widetilde{{\mathcal H}}},
\qquad {\mathcal K}(t,r)={\mathcal K}_{i,j}(t,r)e^{-\lambda(t-r)},
\end{align}
where ${\mathcal K}_{i,j}(t,r)=0$ for $r>t$.

\smallskip
By \eqref{4.4} and \eqref{4.7}, the operator ${\mathcal A}: L^\infty(\mathbb{R};{\mathcal H})\to L^\infty(\mathbb{R})$ is well-defined and \eqref{2.3} holds. Let $k\in\mathbb{N}$ and $f_m\in L^\infty(\mathbb{R};{\mathcal H})$ for $m=1,\ldots,k$. Then
\begin{align}\label{4.24}
\Big|{\mathcal A}\Big(\sum_{m=1}^k f_m(t)\Big)\Big|
=\Big\|{\mathcal A}_{i,j}^\lambda \Big(\sum_{m=1}^k f_m(t,\cdot)\Big)\Big\|_{\widetilde{{\mathcal H}}}
\le \sum_{m=1}^k |{\mathcal A} f_m(t)| \quad \text{a.e.},
\end{align}
so condition (i) of Lemma \ref{lem2.4} holds.

\smallskip
By \eqref{4.22}--\eqref{4.23}, for every $h\in {\mathcal H}$,
\begin{align}\label{4.25}
\begin{aligned}
&|\partial_r{\mathcal K}(t,r)h(x)|
 \\ &\le\bigg|\int_{{\mathbb R}^n}\partial_rK_{i,j}(r,t,y)[h(x-y)-h(x)]dy+
\lambda\int_{{\mathbb R}^n}K_{i,j}(r,t,y)[h(x-y)-h(x)]dy\bigg|e^{-\lambda(t-r)}\\
&\le C[h]_{\alpha,\rho}\!\int_{{\mathbb R}^n}\!(t-r)^{-\frac{n+2}{2}}\big[(t-r)^{-1}+1\big]
e^{-\frac{|y|^2}{4\Gamma(t-r)}}\big(|y|^\alpha\rho(|y|)\mathbf{1}_{\{|y|\le 1\}}+|y|\mathbf{1}_{\{|y|>1\}}\big)dye^{-\lambda(t-r)} \\
&\le C[h]_{\alpha,\rho}\Big[(t-r)^{-\frac{1}{2}}+(t-r)^{\frac{\alpha}{4}-2}\Big] e^{-\lambda(t-r)},
\end{aligned}
\end{align}
where in the penultimate line we used \eqref{2.10}.

\smallskip
For $x\neq y$ with $|x-y|\le 1/3$,
\begin{align}\label{4.26}
\begin{aligned}
\partial_r{\mathcal K}(t,r)h(x)-\partial_r{\mathcal K}(t,r)h(y)
&=\bigg[\int_{{\mathbb R}^n}\partial_rK_{i,j}(r,t,z)[h(x-z)-h(y-z)]dz\\
&\qquad +
\lambda\int_{{\mathbb R}^n}K_{i,j}(r,t,z)[h(x-z)-h(y-z)]dz\bigg]e^{-\lambda(t-r)} \\
&=: [I_1(t,r,x,y)+I_2(t,r,x,y)]e^{-\lambda(t-r)}.
\end{aligned}
\end{align}
Decomposing $I_1$ analogously to \eqref{4.8},
\begin{align}\label{4.27}
\begin{aligned}
I_1(t,r,x,y)
&=\int_{|x-z|\le 2|x-y|}\!\partial^3_{r,x_i,x_j}K(r,t,x-z)[h(z)-h(x)]dz \\
&\quad -\int_{|x-z|\le 2|x-y|}\!\partial^3_{r,y_i,y_j}K(r,t,y-z)[h(z)-h(y)]dz \\
&\quad +\int_{|x-z|> 2|x-y|}\!\partial^3_{r,y_i,y_j}K(r,t,y-z)[h(y)-h(x)]dz \\
&\quad +\int_{2|x-y|<|x-z|\le 1}\!\!\!\big[\partial^3_{r,x_i,x_j}K(r,t,x-z)-\partial^3_{r,y_i,y_j}K(r,t,y-z)\big][h(z)-h(x)]dz \\
&\quad +\int_{|x-z|>1}\!\big[\partial^3_{r,x_i,x_j}K(r,t,x-z)-\partial^3_{r,y_i,y_j}K(r,t,y-z)\big][h(z)-h(x)]dz \\
&=:I_{11}(t,r,x,y)+I_{12}(t,r,x,y)+I_{13}(t,r,x,y)+I_{14}(t,r,x,y)+I_{15}(t,r,x,y).
\end{aligned}
\end{align}
Using that, for every $m\in\mathbb{N}$,
\[
|\partial_r\partial^m_{x_{i_1}\cdots x_{i_m}}K(r,t,x)|\le C(t-r)^{-1-\frac{n+m}{2}}e^{-\frac{|x|^2}{4\Gamma(t-r)}},
\]
and repeating the estimates for $J_1$--$J_5$, we obtain
\begin{align}\label{4.28}
\begin{cases}
|I_{11}(t,r,x,y)|
\le C[h]_{\alpha,\rho}(t-r)^{-2} \displaystyle\int_0^{2|x-y|} \frac{\rho(r)}{r^{1-\alpha}}dr, \\ |I_{12}(t,r,x,y)| \le C[h]_{\alpha,\rho}(t-r)^{-2} \displaystyle\int_0^{3|x-y|} \frac{\rho(r)}{r^{1-\alpha}}dr,\\
|I_{13}(t,r,x,y)| \le C[h]_{\alpha,\rho}(t-r)^{-2}|x-y|^\alpha\rho(|x-y|),\\
|I_{14}(t,r,x,y)| \le C[h]_{\alpha,\rho}(t-r)^{-2}|x-y|\displaystyle\int_{2|x-y|}^1\frac{\rho(r)}{r^{2-\alpha}}dr,\\
|I_{15}(t,r,x,y)| \le C[h]_{\alpha,\rho}(t-r)^{-\frac{3+\beta}{2}}|x-y|^\beta,
\end{cases}
\end{align}
for some $\beta\in(\alpha,1)$. Hence, using \eqref{4.14}--\eqref{4.20} and \eqref{4.28}
\begin{align}\label{4.29}
|I_1(t,r,x,y)|
\le C[h]_{\alpha,\rho}\Big[(t-r)^{-2}+(t-r)^{-\frac{3+\beta}{2}}\Big] |x-y|^\alpha\rho(|x-y|),
\quad |x-y|\le \tfrac{1}{3}.
\end{align}
Similarly,
\begin{align}\label{4.30}
|I_2(t,r,x,y)|
\le C[h]_{\alpha,\rho}\Big[(t-r)^{-1}+(t-r)^{-\frac{1+\beta}{2}}\Big] |x-y|^\alpha\rho(|x-y|),
\quad |x-y|\le \tfrac{1}{3}.
\end{align}
Equations \eqref{4.25}, \eqref{4.26}, \eqref{4.29} and \eqref{4.30} together lead to
\[
\|\partial_r{\mathcal K}(t,r)h\|_{\widetilde{{\mathcal H}}}
\le C\|h\|_{{\mathcal H}}\Big[(t-r)^{-\frac{1}{2}}+(t-r)^{-2}\Big] e^{-\lambda(t-r)},
\]
and therefore
\begin{align}\label{4.31}
\|\partial_r{\mathcal K}(t,r)\|_{L({\mathcal H};\widetilde{{\mathcal H}})}
\le C\Big[(t-r)^{-\frac{1}{2}}+(t-r)^{-2}\Big] e^{-\lambda(t-r)}=:C\phi(t-r).
\end{align}
Similarly,
\begin{align}\label{4.32}
\|{\mathcal K}(t,r)\|_{L({\mathcal H};\widetilde{{\mathcal H}})}
\le C\Big[(t-r)^{-\frac{1}{2}}+(t-r)^{-1}\Big] e^{-\lambda(t-r)}.
\end{align}
By \eqref{4.21}, \eqref{4.23} and \eqref{4.32}, conditions (i)--(ii) of Definition \ref{def2.2} hold.

\smallskip
Noting that, for every $\iota>0$,
\begin{align*}
\iota\int_\iota^{\infty} \phi(\tau)d\tau
=\iota\int_\iota^{\infty}\big[\tau^{-1/2}+\tau^{-2}\big] e^{-\lambda \tau}d\tau\le C<\infty,
\end{align*}
Lemma \ref{lem2.3} yields that ${\mathcal K}$ is a Calder\'{o}n--Zygmund kernel relative to the filtration of dyadic cubes. Moreover, for $f\in {\mathcal C}_0^\infty(\mathbb{R};{\mathcal H})$,
\[
|{\mathcal A} f(t)|
=\bigg\|\int_{\mathbb{R}}e^{-\lambda(t-r)}{\mathcal K}_{i,j}(t,r)f(r,\cdot)dr\bigg\|_{\widetilde{{\mathcal H}}}
= \bigg\|\int_{\mathbb{R}}{\mathcal K}(t,r)f(r)dr\bigg\|_{\widetilde{{\mathcal H}}},
\]
so condition (ii) in Lemma \ref{lem2.4} holds as well.

\smallskip
Let $f,f_1,f_2,\ldots\in L^\infty(\mathbb{R};{\mathcal H})$ vanish outside the same ball. Assume $\|f_k\|_{L^\infty(\mathbb{R};{\mathcal H})}$ are uniformly bounded and $\|f(t)-f_k(t)\|_{{\mathcal H}}\to 0$ for a.e. $t\in\mathbb{R}$. By Fatou's lemma, \eqref{2.4} holds; hence condition (iii) of Lemma \ref{lem2.4} holds. Therefore, by Lemma \ref{lem2.4}, ${\mathcal A}$ is of weak type $(1,1)$ and strong type $(p,p)$ for every $p\in (1,\infty)$ on smooth compactly supported functions. Since ${\mathcal C}_0^\infty(\mathbb{R};{\mathcal H})$ is dense in $L^p(\mathbb{R};{\mathcal H})$, there exists $C>0$ such that
\[
\gamma\big|\{t:\ \|\nabla^2 G_\lambda f(t,\cdot)\|_{\widetilde{{\mathcal H}}}>\gamma\}\big|
\le C\int_{\mathbb{R}}\|f(t,\cdot)\|_{{\mathcal H}}dt, \quad \forall\,\gamma>0,\ \ {\text if} \ \ p=1,
\]
and
\[
\bigg(\int_{\mathbb{R}}\|\nabla^2G_\lambda f(t,\cdot)\|_{\widetilde{{\mathcal H}}}^pdt\bigg)^{\frac{1}{p}}
\le C\bigg(\int_{\mathbb{R}}\|f(t,\cdot)\|_{{\mathcal H}}^pdt\bigg)^{\frac{1}{p}},
\quad {\text if} \ \ p\in (1,\infty),
\]
for all $f\in L^p(\mathbb{R};{\mathcal H})$. Hence \eqref{4.5}--\eqref{4.6} hold. By \eqref{4.1} and \eqref{4.3}--\eqref{4.6}, we conclude \eqref{3.4} and \eqref{3.5} for $\vartheta=l$. \qed

\section{Proof of Theorem \ref{the3.4}}\label{sec5}
\setcounter{equation}{0}
The proof for $\vartheta=b$ is similar to that for
$\vartheta=l$, so we only present the latter case.
The case $p=\infty$ follows directly from
Theorem~\ref{the3.1} on the finite time interval.
We therefore assume $1<p<\infty$ below. Moreover, the cases $2/p-1<\alpha<2/p$ and $\alpha=2/p-1$ are analogous to $\alpha>2/p$ and $\alpha=2/p$, respectively. Hence it suffices to prove the case $\alpha>2/p$ for $\varsigma\in\{d,s,c,w\}$ and the case $\alpha=2/p$ for $\varsigma=d$. When $\varsigma\in\{d,s,c,w\}$, the arguments for $\varsigma\in\{d,c,w\}$ are similar to those for $\varsigma=s$, so we only treat $\varsigma=s$ when $\alpha>2/p$.

\medskip
(i) By the representation \eqref{3.6}, we have
\begin{align}\label{5.1}
\begin{aligned}
\sup_{x\in {\mathbb R}^n}\frac{|u(t,x)|}{1+|x|}
&\le C\sup_{x\in {\mathbb R}^n}\frac{1}{1+|x|}\int_0^t\int_{{\mathbb R}^n}K(r,t,y)[1+|x|+|y|]f_1(r)e^{-\lambda(t-r)}dydr\\
&\le C\int_0^t f_1(r)\big[1+(t-r)^{\frac{1}{2}}\big]dr
\le C\|f_1\|_{L^p([0,T])},
\end{aligned}
\end{align}
where $f_1$ is as in \eqref{4.1}.

\smallskip
For $1\le i\le n$, by \eqref{3.6},
\begin{align}\label{5.2}
\begin{aligned}
|\partial_{x_i}u(t,x)|
&\le C\int_0^t\int_{{\mathbb R}^n}|\partial_{y_i}K(r,t,y)|f_2(r)\big(\mathbf{1}_{\{|y|\le 1\}}|y|^{\alpha}\rho(|y|)+|y|\mathbf{1}_{\{|y|>1\}}\big)e^{-\lambda(t-r)}dydr\\
&\le C\int_0^t\int_{{\mathbb R}^n}(t-r)^{-\frac n2}\exp\Big(-\frac{|y|^2}{C\Gamma(t-r)}\Big)f_2(r)\big[1+|y|^{\alpha-1-\varepsilon}\big]e^{-\lambda(t-r)}dydr\\
&\le C\|f_2\|_{L^p([0,T])}\bigg(1+\Big[\int_0^T r^{\frac{(\alpha-1-\varepsilon)p}{2(p-1)}}dr\Big]^{\frac{p-1}{p}}\bigg)<\infty,
\end{aligned}
\end{align}
where $f_2$ is as in \eqref{4.2}, and in the second line
we used \eqref{2.10} with
$0<\varepsilon<\min\{\alpha,1+\alpha-2/p\}$.

\smallskip
Let $1\le i,j\le n$. Choosing $2\varepsilon=\alpha-2/p$ in \eqref{5.2} yields
\begin{align}\label{5.3}
\begin{aligned}
|\partial^2_{x_i,x_j}u(t,x)|
&\le C\int_0^t f_2(r)\Big[(t-r)^{-\frac{1}{2}}+(t-r)^{\frac{\alpha-\varepsilon}{2}-1}\Big]e^{-\lambda(t-r)}dr \\
&\le C\|f_2\|_{L^p([0,T])}\bigg(\Big[\int_0^T r^{-\frac{p}{2(p-1)}}dr\Big]^{\frac{p-1}{p}}
+\Big[\int_0^T r^{\frac{(\alpha-4)p+2}{4(p-1)}}dr\Big]^{\!\frac{p-1}{p}}\bigg)<\infty,
\end{aligned}
\end{align}
where the last inequality uses $\alpha>2/p$.

\smallskip
Combining \eqref{5.1}--\eqref{5.3}, it remains to prove that for every $1\le i,j\le n$ and every $x,y\in\mathbb{R}^n$ with $|x-y|\le 1/3$, there exists $C>0$ such that
\begin{align}\label{5.4}
|\partial^2_{x_i,x_j}u(t,x)-\partial^2_{y_i,y_j}u(t,y)|
\le C|x-y|^{\alpha-\frac{2}{p}}\rho(|x-y|), \quad t\in[0,T].
\end{align}

\smallskip
Let $J_1,\ldots,J_5$ be defined as in \eqref{4.8} with $-\infty$ replaced by $0$. By H\"older's inequality,
\begin{align}\label{5.5}
\begin{aligned}
|J_1(t,x,y)|
&\le C[f]_{p,\alpha,\rho}\int_{|z|\le 2|x-y|}\!|z|^\alpha\rho(|z|)\bigg[\int_0^{\infty} r^{-\frac{(n+2)p}{2(p-1)}}e^{-\frac{p|z|^2}{4\Gamma(p-1)r}}dr\bigg]^{\frac{p-1}{p}}dz \\
&\le C\int_{|z|\le 2|x-y|}\frac{\rho(|z|)}{|z|^{n+\frac{2}{p}-\alpha}}dz
\le C\int_0^{2|x-y|}\frac{\rho(r)}{r^{1+\frac{2}{p}-\alpha}}dr,
\end{aligned}
\end{align}
where $[f]_{p,\alpha,\rho}^p=\int_0^T [f(r,\cdot)]_{\alpha,\rho}^pdr$.

\smallskip
Similarly,
\begin{align}\label{5.6}
|J_2(t,x,y)|\le C\int_0^{3|x-y|}\frac{\rho(r)}{r^{1+\frac{2}{p}-\alpha}}dr.
\end{align}
For $J_3$,
\begin{align}\label{5.7}
\begin{aligned}
|J_3(t,x,y)|
&\le C[f]_{p,\alpha,\rho}|x-y|^\alpha\rho(|x-y|)|x-y|^n
\bigg[\int_0^{\infty} r^{-\frac{(n+2)p}{2(p-1)}}e^{-\frac{p|x-y|^2}{2\Gamma(p-1)r}}dr\bigg]^{\!\frac{p-1}{p}} \\
&\le C|x-y|^{\alpha-\frac{2}{p}}\rho(|x-y|).
\end{aligned}
\end{align}
For $J_4$,
\begin{align}\label{5.8}
\begin{aligned}
|J_4(t,x,y)|
&\le C[f]_{p,\alpha,\rho}|x-y|
\int_{2|x-y|<|z|\le 1}\!|z|^\alpha\rho(|z|)\bigg[\int_0^{\infty} r^{-\frac{(n+3)p}{2(p-1)}}e^{-\frac{p|z|^2}{16\Gamma(p-1)r}}dr\bigg]^{\!\frac{p-1}{p}}dz\\
&\le C|x-y|\int_{2|x-y|}^1\frac{\rho(r)}{r^{2+\frac{2}{p}-\alpha}}dr.
\end{aligned}
\end{align}
Finally,
\begin{align}\label{5.9}
\begin{aligned}
|J_5(t,x,y)|
&\le C[f]_{p,\alpha,\rho}|x-y|^\beta
\int_{|z|>1}\!|z|\bigg[\int_0^{\infty} r^{-\frac{(n+2+\beta)p}{2(p-1)}}e^{-\frac{p|z|^2}{16\Gamma(p-1)r}}dr\bigg]^{\!\frac{p-1}{p}}dz\\
&\le C|x-y|^\beta\int_{|z|>1}|z|^{-n-\beta-\frac{2}{p}+1}dz
\le C|x-y|^\beta,
\end{aligned}
\end{align}
for some $\beta\in(\alpha,1)$ with $\beta+2/p>1$.

\smallskip
Combining \eqref{5.5}--\eqref{5.9} and \eqref{2.10}, for every $t\in[0,T]$ it follows that
\begin{align}\label{5.10}
\begin{aligned}
&|\partial^2_{x_i,x_j}u(t,x)-\partial^2_{y_i,y_j}u(t,y)|\\
&\le C\bigg[\int_0^{3|x-y|}\frac{\rho(r)}{r^{1+\frac{2}{p}-\alpha}}dr
+|x-y|^{\alpha-\frac{2}{p}}\rho(|x-y|)
+|x-y|\int_{2|x-y|}^1\frac{\rho(r)}{r^{2+\frac{2}{p}-\alpha}}dr\bigg].
\end{aligned}
\end{align}
Adapting the computations for \eqref{4.14}--\eqref{4.20} to \eqref{5.10}, we conclude \eqref{5.4}.

\medskip
(ii) From \eqref{5.5}--\eqref{5.10} and \eqref{4.15} (with $\alpha=0$), it follows that, for $|x-y|\le 1/3$,
\begin{align}\label{5.11}
|\nabla^2u(t,x)-\nabla^2u(t,y)|
\le C\|f\|_{L^p([0,T];{\mathcal C}_{l,d}^{\alpha,\rho}({\mathbb R}^n))}\hat{\rho}(|x-y|),
\quad t\in[0,T].
\end{align}

On the other hand, by \eqref{5.1} and \eqref{5.2}, we have $(1+|\cdot|)^{-1}u(t,\cdot),\,\partial_{x_i}u\in L^\infty([0,T]\times{\mathbb R}^n)$. It remains to show $\partial^2_{x_i,x_j}u\in L^\infty([0,T]\times{\mathbb R}^n)$ for $1\le i,j\le n$. Indeed,
\begin{align}\label{5.12}
\begin{aligned}
\big|\partial^2_{x_i,x_j}u(t,x)\big|&=\bigg|\int_0^t\int_{{\mathbb R}^n}e^{-\lambda(t-r)}\partial^2_{x_i,x_j}K(r,t,x-y)[f(r,y)-f(r,x)]dydr\bigg| \\
&\le C\|f\|_{L^p([0,T];{\mathcal C}_{l,d}^{\alpha,\rho}({\mathbb R}^n))}
\int_{{\mathbb R}^n}\bigg[\int_0^T r^{-\frac{p(n+2)}{2(p-1)}}e^{-\frac{p|y|^2}{8\Gamma(p-1)r}}dr\bigg]^{\frac{p-1}{p}}
\\ &\qquad \times
\big(|y|^\alpha\rho(|y|)\mathbf{1}_{\{|y|\le 1\}}+|y|\mathbf{1}_{\{|y|>1\}}\big)e^{-\frac{|y|^2}{8\Gamma T}}dy \\
&\le C\|f\|_{L^p([0,T];{\mathcal C}_{l,d}^{\alpha,\rho}({\mathbb R}^n))}
\bigg[\int_0^1\frac{\rho(r)}{r}dr+\int_1^\infty r^{-\alpha}e^{-\frac{r^2}{8\Gamma T}}dr\bigg]
\\
&\le C\|f\|_{L^p([0,T];{\mathcal C}_{l,d}^{\alpha,\rho}({\mathbb R}^n))},
\end{aligned}
\end{align}
which completes the proof. \qed

\section{Proof of Theorem \ref{the3.6}}\label{sec6}
\setcounter{equation}{0}
We first establish the unique solvability and the maximal regularity
estimate \eqref{3.13} for $\vartheta\in\{l,b\}$ and
$\varsigma\in\{d_k,s,c,w\}$.
The additional time-uniform estimates are proved later
by applying Theorem~\ref{the3.4} and Remark~\ref{rem3.5}
to a spatially frozen equation.
We only treat the case $\vartheta=l$ and $\alpha>2/p$, since the case $\vartheta=b$ or $\alpha\le 2/p$ are analogous and simpler. Let ${\mathcal X}$ be given by \eqref{3.12}. Suppose $u\in {\mathcal X}$
solves the Cauchy problem \eqref{3.11}. For any $\tilde{\lambda}\in\mathbb{R}$, define
$\tilde{u}(t,x):=u(t,x)e^{(\lambda-\tilde{\lambda})t}$.
Then $\tilde{u}$ belongs to the same space and satisfies
\[
\left\{
\begin{aligned}
&\partial_t \tilde{u}(t,x)= \frac12\sum_{i,j=1}^n a_{i,j}(t,x)\partial^2_{x_i,x_j}\tilde{u}(t,x)
+ g(t,x)\cdot\nabla \tilde{u}(t,x) \\
&\qquad\qquad \ \ - \tilde{\lambda}\tilde{u}(t,x) + \tilde{f}(t,x), \quad (t,x)\in(0,T]\times\mathbb{R}^n, \\
&\tilde{u}(0,x)=0, \quad x\in\mathbb{R}^n,
\end{aligned}
\right.
\]
with $\tilde{f}(t,x)=f(t,x)e^{(\lambda-\tilde{\lambda})t}$, and the converse also holds. Hence it suffices to establish well-posedness of \eqref{3.11} in $\mathcal X$ for sufficiently large $\lambda$. In this case $u$ admits the representation
\begin{align}\label{6.1}
u(t,x)=\int_0^t\int_{\mathbb{R}^n} e^{-\lambda(t-r)}\hat{K}(r,t,y,x)
\big[g(r,y)\cdot\nabla u(r,y)+f(r,y)\big]dydr,
\end{align}
where $\hat K(r,t,y,x)$ is the fundamental solution in
$(t,x)$ with initial pole $(r,y)$ for the operator
\[
\partial_t-\tfrac12\sum_{i,j=1}^n a_{i,j}(t,x)\partial^2_{x_i,x_j}.
\]
By \eqref{2.10}, the coefficients are uniformly H\"older continuous
in space with any fixed exponent $\gamma\in(0,\alpha)$.
Therefore, \cite[Theorems~4.6 and~4.7]{PP20} apply.
For $t>r$, the kernel is twice continuously differentiable in $x$
and satisfies the following Gaussian estimates:
\begin{align}\label{6.2}
\left\{
\begin{aligned}
&\frac{C_1}{(t-r)^{\frac n2}}
 \exp\Big(-\frac{C_2|x-y|^2}{t-r}\Big)
 \le \hat K(r,t,y,x)
 \le \frac{C_3}{(t-r)^{\frac n2}}
 \exp\Big(-\frac{C_4|x-y|^2}{t-r}\Big),\\
&\big|\nabla_x^m\hat K(r,t,y,x)\big|
 \le \frac{C_5}{(t-r)^{\frac{n+m}{2}}}
 \exp\Big(-\frac{C_6|x-y|^2}{t-r}\Big),
 \qquad m=1,2,
\end{aligned}
\right.
\end{align}
where $0\le r<t\le T$, $x,y\in\mathbb R^n$,
and $C_i$, $i=1,\ldots,6$, are positive constants.

\medskip
For $\tau\in[0,1]$, consider
\begin{align}\label{6.3}
\left\{
\begin{aligned}
&\partial_t u(t,x)= \frac12\sum_{i,j=1}^n a_{i,j}(t,x)\partial^2_{x_i,x_j}u(t,x)
+ \tau g(t,x)\cdot\nabla u(t,x) \\
&\qquad\qquad \ \ - \lambda u(t,x) + f(t,x), \quad (t,x)\in(0,T]\times\mathbb{R}^n,\\
&u(0,x)=0,\quad x\in\mathbb{R}^n.
\end{aligned}
\right.
\end{align}

When $\tau=0$, \eqref{6.3} reduces to
\begin{align}\label{6.4}
\left\{
\begin{aligned}
&\partial_t u(t,x)= \frac12\sum_{i,j=1}^n a_{i,j}(t,x)\partial^2_{x_i,x_j}u(t,x)
 - \lambda u(t,x) + f(t,x), \quad (t,x)\in(0,T]\times\mathbb{R}^n,\\
&u(0,x)=0,\quad x\in\mathbb{R}^n.
\end{aligned}
\right.
\end{align}
For the drift-free equation \eqref{6.4}, we establish
the a priori estimates and unique solvability for every
$1<p\le\infty$. This full range will also be used in (ii).

\smallskip
By Remark~\ref{rem2.11}, for $x\ne y$,
\begin{align}\label{6.5}
 |f(r,x)-f(r,y)|\le C[f(r,\cdot)]_{\alpha,\rho}\bigl(
 |x-y|^\alpha\rho(|x-y|)
       \mathbf1_{\{0<|x-y|\le1\}}
 +|x-y|\mathbf1_{\{|x-y|>1\}}
 \bigr).
\end{align}
We use the same decomposition for the coefficient increments
below.

\smallskip
Fix $0<\varepsilon<\min\{\alpha,1-\alpha\}$.
The representation \eqref{2.7} and the limits \eqref{2.8} in Lemma~\ref{lem2.8} give
\begin{align}\label{6.6}
 \frac{\rho(a)}{\rho(b)}
 \le C_\varepsilon
 \max\left\{
       \left(\frac ab\right)^{-\varepsilon},
       \left(\frac ab\right)^\varepsilon
 \right\},
 \qquad 0<a,b\le1.
\end{align}
For $0<v\le1$ and $C_0>0$, changing variables $z=\sqrt vw$ gives
\begin{align*}
\int_{|z|\le1}v^{-\frac n2}e^{-\frac{C_0|z|^2}{v}}
       |z|^\alpha\rho(|z|)dz
 &\le
 C v^{\frac\alpha2}\rho(\sqrt v)
 \int_{\mathbb R^n}e^{-C_0|w|^2}
       \bigl(|w|^{\alpha-\varepsilon}
             +|w|^{\alpha+\varepsilon}\bigr)dw\le C v^{\frac\alpha2}\rho(\sqrt v).
\end{align*}
Moreover, by \eqref{2.10},
\[
 \int_{|z|>1}v^{-\frac n2}e^{-\frac{C_0|z|^2}{v}}|z|dz
 \le C\sqrt v
 \le C v^{\frac\alpha2}\rho(\sqrt v).
\]
For $v>1$, boundedness of $r^\alpha\rho(r)$ on $(0,1]$
and the Gaussian first-moment estimate give
\begin{align}\label{6.7}
 &\int_{\mathbb R^n}v^{-\frac n2}e^{-\frac{C_0|z|^2}{v}}
 \bigl(
 |z|^\alpha\rho(|z|)\mathbf1_{\{0<|z|\le1\}}
 +|z|\mathbf1_{\{|z|>1\}}
 \bigr)dz\le
 C\begin{cases}
 v^{\frac\alpha2}\rho(\sqrt v),&0<v\le1,\\
 \sqrt v,&v>1.
 \end{cases}
\end{align}

Since $\int_{\mathbb R^n}\hat K(r,t,y,x)dy=1$,
the derivative bounds in \eqref{6.2} imply
\begin{align}\label{6.8}
 \int_{\mathbb R^n}\nabla_x^m\hat K(r,t,y,x)dy=0,
 \qquad m=1,2.
\end{align}
Using this cancellation \eqref{6.8} and \eqref{6.2}, and arguing as in
\eqref{4.1}--\eqref{4.4}, we obtain from \eqref{6.1} with $g=0$
\begin{align*}
\left\{
\begin{aligned}
 &\frac{|u(t,x)|}{1+|x|}\le C\int_0^t e^{-\lambda(t-r)}
       \big(1+\sqrt{t-r}\big)
       \|(1+|\cdot|)^{-1}f(r,\cdot)\|_0dr,\\
 &\|\nabla^m u(t,\cdot)\|_0
 \le C_T\int_0^t e^{-\lambda(t-r)}
 \bigl[1+(t-r)^{\frac\alpha4-\frac m2}\bigr]
 [f(r,\cdot)]_{\alpha,\rho}dr,
 \quad m=1,2.
\end{aligned}
\right.
\end{align*}
These time kernels are integrable. By Young's
inequality, we obtain, uniformly for $\lambda\geq1$,
\begin{align}\label{6.9}
\begin{aligned}
 &\big\|\|(1+|\cdot|)^{-1}u(t,\cdot)\|_0\big\|_{L^p(0,T)}
 +\big\|\|\nabla u(t,\cdot)\|_0\big\|_{L^p(0,T)}\\
 &+
 \big\|\|\nabla^2u(t,\cdot)\|_0\big\|_{L^p(0,T)}
 \le C\|f\|_{L^p([0,T];
       \mathcal C_{l,\varsigma}^{\alpha,\rho}(\mathbb R^n))}.
\end{aligned}
\end{align}

For fixed $x_0\in\mathbb R^n$, let $K$ be the Gaussian kernel
\eqref{3.2} associated with the matrix $a(\cdot,x_0)$. Then
\begin{align}\label{6.10}
u(t,x)=\int_0^t\int_{\mathbb R^n}
 e^{-\lambda(t-r)}K(r,t,x-z)
 \bigg[f(r,z)
 +\frac12\sum_{i,j=1}^n
 \bigl(a_{i,j}(r,z)-a_{i,j}(r,x_0)\bigr)
 \partial^2_{z_i,z_j}u(r,z)\bigg]dzdr.
\end{align}
Fix $x_1$ with $|x_1-x_0|\in(0,1]$. For $0<|z-x_0|\le1$,
\begin{align}\label{6.11}
 &\big|
 \bigl(a_{i,j}(r,z)-a_{i,j}(r,x_0)\bigr)
 \partial^2_{z_i,z_j}u(r,z)
 \big|\le
 C|z-x_0|^\alpha\rho(|z-x_0|)
 \|\nabla^2u(r,\cdot)\|_0.
\end{align}
For $0<|z-x_1|\le1$,
\begin{align}\label{6.12}
\begin{aligned}
&\big|\big(a_{i,j}(r,z)-a_{i,j}(r,x_0)\big)
       \partial^2_{z_i,z_j}u(r,z)
-\big(a_{i,j}(r,x_1)-a_{i,j}(r,x_0)\big)
       \partial^2_{x_i,x_j}u(r,x_1)\big|\\
&\le
|a_{i,j}(r,z)-a_{i,j}(r,x_1)|
\|\nabla^2u(r,\cdot)\|_0+
|a_{i,j}(r,x_1)-a_{i,j}(r,x_0)|
\big|\partial^2_{z_i,z_j}u(r,z)
     -\partial^2_{x_i,x_j}u(r,x_1)\big|\\
&\le
C|z-x_1|^\alpha\rho(|z-x_1|)
\big(\|\nabla^2u(r,\cdot)\|_0
+|x_1-x_0|^\alpha\rho(|x_1-x_0|)[\nabla^2u(r,\cdot)]_{\alpha,\rho}\big),
\end{aligned}
\end{align}
where the seminorm on the right hand side is finite for almost every $r$
because the present calculation concerns an a priori estimate
for $u\in\mathcal X$.
For $|z-x_0|>1$ in \eqref{6.11}, or $|z-x_1|>1$
in \eqref{6.12}, the corresponding left-hand side is bounded
by $C\|\nabla^2u(r,\cdot)\|_0$, since the coefficients
$a_{i,j}$ are bounded.

\smallskip
For $t-r\geq h^2$, the mean value theorem and the
third-derivative bound for the explicit frozen Gaussian give
\begin{align}\label{6.13}
|\nabla^2K(r,t,x_1-z)-\nabla^2K(r,t,x_0-z)|\le
 C|x_1-x_0|(t-r)^{-\frac{n+3}{2}}
 \exp\biggl(-\frac{\hat C|x_0-z|^2}{t-r}\biggr),
\end{align}
where $\hat C>0$ is a constant. Split the time integral in \eqref{6.10} into the regions
$t-r<|x_1-x_0|^2$, $|x_1-x_0|^2\le t-r\le1$, and $t-r>1$.
Since
\[
 \int_{\mathbb R^n}\nabla^2K(r,t,x-z)dz=0,
\]
in the region $t-r<|x_1-x_0|^2$, we estimate the two Hessians
separately, using the increments
\[
 f(r,z)-f(r,x_0)
 \quad\text{at }x_0,
 \qquad
 f(r,z)-f(r,x_1)
 \quad\text{at }x_1,
\]
and treating the coefficient terms by
\eqref{6.11} and \eqref{6.12}, respectively.
In the regions $t-r\geq |x_1-x_0|^2$, we combine the two Hessians
and use \eqref{6.13}, with $f(r,z)$ replaced by
$f(r,z)-f(r,x_0)$. The coefficient perturbation already
vanishes at $z=x_0$.
For the term involving $[\nabla^2u]_{\alpha,\rho}$,
we also use $s^\alpha\rho(s)\le Cs^{\alpha/2}$
for $0<s\le1$,
which follows from \eqref{2.10}. Combining \eqref{6.7} and \eqref{6.10}--\eqref{6.13},
we obtain
\begin{align}\label{6.14}
\begin{aligned}
&|\nabla^2u(t,x_1)-\nabla^2u(t,x_0)|
\\ &\le
 C\int_{(t-|x_1-x_0|^2)_+}^t
 (t-r)^{\frac\alpha2-1}\rho(\sqrt{t-r})
 \bigl(\|f(r,\cdot)\|_{l,\alpha,\rho}
       +\|\nabla^2u(r,\cdot)\|_0\bigr)dr\\
 &\quad+
 C|x_1-x_0|\int_{(t-1)_+}^{(t-|x_1-x_0|^2)_+}
 (t-r)^{\frac{\alpha-3}{2}}\rho(\sqrt{t-r})
 \bigl(\|f(r,\cdot)\|_{l,\alpha,\rho}
       +\|\nabla^2u(r,\cdot)\|_0\bigr)dr\\
 &\quad+
 C_T |x_1-x_0|\int_0^{(t-1)_+}
 \bigl(\|f(r,\cdot)\|_{l,\alpha,\rho}
       +\|\nabla^2u(r,\cdot)\|_0\bigr)dr\\
 &\quad+
 C|x_1-x_0|^\alpha\rho(|x_1-x_0|)\int_0^t
 e^{-\lambda(t-r)}(t-r)^{\frac\alpha4-1}
 [\nabla^2u(r,\cdot)]_{\alpha,\rho}dr.
\end{aligned}
\end{align}

For $0<v\le |x_1-x_0|^2$, \eqref{6.6} gives
\begin{align}\label{6.15}
 \frac{v^{\alpha/2}\rho(\sqrt v)}{|x_1-x_0|^\alpha\rho(|x_1-x_0|)}
 \le C_\varepsilon
 \left(\frac{v}{|x_1-x_0|^2}\right)^{\frac{\alpha-\varepsilon}{2}}.
\end{align}
For $|x_1-x_0|^2\le v\le1$, the same estimate gives
\begin{align}\label{6.16}
 \frac{v^{\alpha/2}\rho(\sqrt v)}{|x_1-x_0|^\alpha\rho(|x_1-x_0|)}
 \le C_\varepsilon
 \left(\frac{v}{|x_1-x_0|^2}\right)^{\frac{\alpha+\varepsilon}{2}}.
\end{align}
Furthermore, \eqref{2.10} implies
\[
 \frac{|x_1-x_0|}{|x_1-x_0|^\alpha\rho(|x_1-x_0|)}\le C,
 \qquad 0<|x_1-x_0|\le1.
\]
After division by $|x_1-x_0|^\alpha\rho(|x_1-x_0|)$, the third term in
\eqref{6.14} is therefore bounded by
\[
 C_T\int_0^t
 \bigl(\|f(r,\cdot)\|_{l,\alpha,\rho}
       +\|\nabla^2u(r,\cdot)\|_0\bigr)dr.
\]
Using dyadic time decompositions for the first two terms,
we obtain
\begin{align}\label{6.17}
\begin{aligned}
\frac{|\nabla^2u(t,x_1)-\nabla^2u(t,x_0)|}
        {|x_1-x_0|^\alpha\rho(|x_1-x_0|)}
&\le
 C_T\bigg(
 \sum_{k=0}^\infty2^{-k(\alpha-\varepsilon)/2}
 +\sum_{k=0}^\infty2^{-k(1-\alpha-\varepsilon)/2}
 \bigg)\\
&\quad\times
 \sup_{\delta>0}\frac1\delta
 \int_{(t-\delta)^+}^t
 \bigl(\|f(r,\cdot)\|_{l,\alpha,\rho}
       +\|\nabla^2u(r,\cdot)\|_0\bigr)dr\\
&\quad
 +C\int_0^t e^{-\lambda(t-r)}(t-r)^{\frac\alpha4-1}
 [\nabla^2u(r,\cdot)]_{\alpha,\rho}dr.
\end{aligned}
\end{align}
Both series converge. Taking the spatial supremum and applying the one-sided
Hardy--Littlewood maximal inequality and Young's inequality,
we obtain
\begin{align}\label{6.18}
\begin{aligned}
 \big\|[\nabla^2u(t,\cdot)]_{\alpha,\rho}\big\|_{L^p(0,T)}
 &\le C\|f\|_{L^p([0,T];
       \mathcal C_{l,\varsigma}^{\alpha,\rho}(\mathbb R^n))}
 +C\big\|\|\nabla^2u(t,\cdot)\|_0\big\|_{L^p(0,T)}\\
 &\quad+
 C\bigg(\int_0^T e^{-\lambda v}v^{\frac\alpha4-1}dv\bigg)
 \big\|[\nabla^2u(t,\cdot)]_{\alpha,\rho}\big\|_{L^p(0,T)}.
\end{aligned}
\end{align}
Since
\[
 \lim_{\lambda\to\infty}
 \int_0^T e^{-\lambda v}v^{\frac\alpha4-1}dv=0,
\]
we choose $\lambda$ sufficiently large to absorb the last term.
Combining \eqref{6.9}, \eqref{6.18} and \eqref{6.4},
we conclude that
\begin{align}\label{6.19}
\begin{aligned}
 \|\partial_tu\|_{L^p([0,T];
       \mathcal C_{l,\varsigma}^{\alpha,\rho}(\mathbb R^n))}
 +\|u\|_{L^p([0,T];
       \mathcal C_{l,\varsigma}^{2+\alpha,\rho}(\mathbb R^n))}\le
 C\|f\|_{L^p([0,T];
       \mathcal C_{l,\varsigma}^{\alpha,\rho}(\mathbb R^n))},
 \qquad 1< p\le\infty.
\end{aligned}
\end{align}
The same estimate holds on every interval $[0,t]\subset[0,T]$,
with a constant independent of $t$.

\smallskip
For the subsequent time-uniform estimates, we fix $\lambda$
large enough and construct a solution of \eqref{6.4} in $\mathcal X$.
Set
\[
 a_{i,j}^{(\theta)}(t,x)
 =(1-\theta)\delta_{i,j}+\theta a_{i,j}(t,x),
 \qquad 0\le\theta\le1,
\]
where $\delta_{i,j}$ is the Kronecker delta.
The ellipticity constants and spatial coefficient norms are
uniform in $\theta$. Thus the preceding a priori argument \eqref{6.19} gives
\begin{align}\label{6.20}
\|\partial_tu\|_{L^p([0,T];
       \mathcal C_{l,\varsigma}^{\alpha,\rho}(\mathbb R^n))}
 +\|u\|_{L^p([0,T];
       \mathcal C_{l,\varsigma}^{2+\alpha,\rho}(\mathbb R^n))}\le
 C\bigg\|
 \partial_tu-\frac12\sum_{i,j=1}^n
 a_{i,j}^{(\theta)}\partial^2_{x_i,x_j}u+\lambda u
 \bigg\|_{L^p([0,T];
       \mathcal C_{l,\varsigma}^{\alpha,\rho}(\mathbb R^n))}
\end{align}
for every $u\in\mathcal X$ with $u(0,x)=0$, with the same
sufficiently large $\lambda$ and a constant independent of
$\theta$.

\smallskip
At $\theta=0$, the unique solvability of \eqref{6.4} follows from
Theorem~\ref{the3.1}.
Suppose now that the unique solvability in $\mathcal X$ has been established at
$\theta\in[0,1)$. Set $u_0=0$ and, for
$0<\delta\le1-\theta$, define
\begin{align}\label{6.21}
\left\{
\begin{aligned}
 &\partial_tu_{k+1}(t,x)
 =\frac12\sum_{i,j=1}^n
 a_{i,j}^{(\theta)}(t,x)\partial^2_{x_i,x_j}u_{k+1}(t,x)
 -\lambda u_{k+1}(t,x)+f(t,x)\\
 &\qquad\qquad\qquad
 +\frac{\delta}{2}\sum_{i,j=1}^n
 (a_{i,j}(t,x)-\delta_{i,j})
 \partial^2_{x_i,x_j}u_k(t,x),\\
 &u_{k+1}(0,x)=0.
\end{aligned}
\right.
\end{align}
The product estimate gives
\begin{align}\label{6.22}
\bigg\|\sum_{i,j=1}^n
 (a_{i,j}-\delta_{i,j})\partial^2_{x_i,x_j}u_k
 \bigg\|_{L^p([0,T];
       \mathcal C_{l,\varsigma}^{\alpha,\rho}(\mathbb R^n))}\le
 C\|u_k\|_{L^p([0,T];
       \mathcal C_{l,\varsigma}^{2+\alpha,\rho}(\mathbb R^n))}.
\end{align}
Starting from $u_0=0$, solvability at $\theta$ and
\eqref{6.22} show inductively that every $u_k$ belongs to
$\mathcal X$. Applying \eqref{6.20} to consecutive differences
in \eqref{6.21} gives
\[
 \|u_{k+1}-u_k\|_{\mathcal X}
 \le C\delta\|u_k-u_{k-1}\|_{\mathcal X},
 \qquad k\geq1.
\]
Choose $\delta>0$ sufficiently small that $C\delta\le1/2$.
Then $u_k$ converges in $\mathcal X$. Since
\[
 a_{i,j}^{(\theta+\delta)}
 =a_{i,j}^{(\theta)}
   +\delta(a_{i,j}-\delta_{i,j}),
\]
passing to the limit in \eqref{6.21} gives a solution with
diffusion matrix $a^{(\theta+\delta)}$. Uniqueness follows
from \eqref{6.20}. The admissible step size is independent
of $\theta$, so finitely many steps give a unique solution
of \eqref{6.4} in $\mathcal X$.

\smallskip
Fix $x_0\in\mathbb R^n$. The solution of \eqref{6.4} also
satisfies the spatially frozen representation \eqref{6.10}.
The product estimate and \eqref{6.19} with $\theta=1$ imply, for every $1<p\le\infty$, that
\begin{align}\label{6.23}
\begin{aligned}
 \bigg\|
 f+\frac12\sum_{i,j=1}^n
 \bigl(a_{i,j}(\cdot,\cdot)-a_{i,j}(\cdot,x_0)\bigr)
       \partial^2_{x_i,x_j}u
 \bigg\|_{L^p([0,T];
       \mathcal C_{l,\varsigma}^{\alpha,\rho}(\mathbb R^n))}\le
 C\|f\|_{L^p([0,T];
       \mathcal C_{l,\varsigma}^{\alpha,\rho}(\mathbb R^n))}.
\end{aligned}
\end{align}
The constants are independent of $x_0$.

\smallskip
The frozen diffusion matrix associated with \eqref{6.10} depends only on time.
Applying Theorem~\ref{the3.4} (i) and \eqref{6.23}, yield
\begin{align}\label{6.24}
\|u\|_{L^\infty([0,T];
       \mathcal C_{l,\varsigma}^{2+\alpha-\frac{2}{p},\rho}(\mathbb R^n))}\le
 C\|f\|_{L^p([0,t];
       \mathcal C_{l,\varsigma}^{\alpha,\rho}(\mathbb R^n))},
\quad 1< p\le \infty.
\end{align}

\smallskip
(i) Assume $q=p\in[2,\infty]$.  Define
\[
{\mathcal H}_1:=L^\infty([0,T];{\mathcal C}_{l,\varsigma}^{2+\alpha-\frac{2}{p},\rho}(\mathbb{R}^n)).
\]
Let $\tau=\tau_0\in (0,1)$ and $v\in \mathcal X\cap\mathcal H_1$. Consider the following Cauchy problem
\begin{align}\label{6.25}
\left\{
\begin{aligned}
&\partial_t u(t,x)= \frac12\sum_{i,j=1}^n a_{i,j}(t,x)\partial^2_{x_i,x_j}u(t,x)
- \lambda u(t,x)  \\
&\qquad\qquad \ \ + \tau_0 g(t,x)\cdot\nabla v(t,x) + f(t,x), \quad (t,x)\in(0,T]\times\mathbb{R}^n,\\
&u(0,x)=0,\quad x\in\mathbb{R}^n.
\end{aligned}
\right.
\end{align}
By the conclusion establishing for the drift-free equation \eqref{6.3}, there is a unique $\mathcal S(v)\in \mathcal X\cap\mathcal H_1$ solving \eqref{6.25}. Moreover, for $v_1,v_2\in \mathcal X\cap\mathcal H_1$, by \eqref{6.20} and \eqref{6.24}, we have
\begin{align}\label{6.26}
\begin{aligned}
 \|\mathcal S(v_1)-\mathcal S(v_2)\|_{\mathcal X\cap\mathcal H_1}
 &\le C\tau_0
 \|g\cdot\nabla(v_1-v_2)\|_{L^p([0,T];
 \mathcal C_{l,\varsigma}^{\alpha,\rho}(\mathbb R^n))}\\
 &\le C\tau_0
 \|g\|_{L^p([0,T];
 \mathcal C_{b,\varsigma}^{\alpha,\rho}(\mathbb R^n))}
 \|v_1-v_2\|_{\mathcal X\cap\mathcal H_1}.
\end{aligned}
\end{align}
Choose $\tau_0$ sufficiently small that
\[
 C\tau_0
 \|g\|_{L^p([0,T];
 \mathcal C_{b,\varsigma}^{\alpha,\rho}(\mathbb R^n))}
 \le\frac12.
\]
Then $\mathcal S$ is a contraction on $\mathcal X\cap\mathcal H_1$.
Its unique fixed point belongs to
$\mathcal X\cap\mathcal H_1$ and solves \eqref{6.3}
at $\tau=\tau_0$. Moreover, \eqref{6.19} and \eqref{6.24} hold true as well.

\smallskip
Suppose that \eqref{6.3} is uniquely solvable in
$\mathcal X\cap\mathcal H_1$ at some $\tau_0\in (0,1)$.
For $0<\delta\le1-\tau_0$ and
$v\in\mathcal X\cap\mathcal H_1$, consider
\begin{align}\label{6.27}
\left\{
\begin{aligned}
 &\partial_tu(t,x)
 =\frac12\sum_{i,j=1}^n
 a_{i,j}(t,x)\partial^2_{x_i,x_j}u(t,x)
 +\tau_0g(t,x)\cdot\nabla u(t,x)\\
 &\qquad\qquad
 -\lambda u(t,x)
 +\delta g(t,x)\cdot\nabla v(t,x)+f(t,x),\\
 &u(0,x)=0.
\end{aligned}
\right.
\end{align}
Then \eqref{6.27} has a unique solution
$\mathcal S(v)\in\mathcal X\cap\mathcal H_1$.
For $v_1,v_2\in\mathcal X\cap\mathcal H_1$,
applying the uniform a priori estimate to the difference
equation, with drift $\tau_0g(t,x)\cdot\nabla u(t,x)$ and forcing term
$\delta g\cdot\nabla(v_1-v_2)$, yields the following
analogue of \eqref{6.26}:
\[
 \|\mathcal S(v_1)-\mathcal S(v_2)\|
       _{\mathcal X\cap\mathcal H_1} \le
 C\delta
 \|g\|_{L^p([0,T];
 \mathcal C_{b,\varsigma}^{\alpha,\rho}(\mathbb R^n))}
 \|v_1-v_2\|_{\mathcal X\cap\mathcal H_1},
\]
where $C$ is independent of $\tau_0$. Choose $\delta_*\in(0,1]$ such that
\[
 C\delta_*
 \|g\|_{L^p([0,T];
 \mathcal C_{b,\varsigma}^{\alpha,\rho}(\mathbb R^n))}
 \le\frac12,
\]
and take $\delta=\min\{\delta_*,1-\tau_0\}$.
Then $\mathcal S$ is a contraction.
Its unique fixed point solves \eqref{6.3}
at $\tau=\tau_0+\delta$. Since $\delta_*$ is independent of $\tau_0$,
starting from $\tau=0$ and repeating this construction
reaches $\tau=1$ after finitely many steps. Moreover, by \eqref{6.19} and \eqref{6.24}, \eqref{3.13} and \eqref{3.8} hold.

\medskip
(ii) Now assume $q=\infty$.
For every $u\in\mathcal X$, the product estimate gives
\[
 \|g\cdot\nabla u\|_{L^p([0,T];
 \mathcal C_{l,\varsigma}^{\alpha,\rho}(\mathbb R^n))}
 \le
 C\|g\|_{L^\infty([0,T];
 \mathcal C_{b,\varsigma}^{\alpha,\rho}(\mathbb R^n))}
 \|u\|_{L^p([0,T];
 \mathcal C_{l,\varsigma}^{2+\alpha,\rho}(\mathbb R^n))}.
\]
Using the drift-free estimates established above, we repeat
the a priori estimates and fixed-point construction in (i),
now in $\mathcal X$, to obtain unique solvability and
\eqref{3.13}. Applying Theorem~\ref{the3.4} (i) and \eqref{6.24}
 to the spatially frozen equation give the desired conclusion. \qed

\section{Proof of Theorem \ref{the1.2}}\label{sec7}\setcounter{equation}{0}
The proof is divided into two parts:
first, we establish the unique strong solvability of a class of Kolmogorov equations;
second, we apply the resulting transform to derive the well-posedness and flow property
for SDE \eqref{1.1}. For simplicity, throughout this section we assume $T\le 1$.

\medskip
\textbf{Part I:} The unique strong solvability for the Kolmogorov equation
\begin{align}\label{7.1}
\left\{
\begin{aligned}{}
&\partial_tU(t,x)=\tfrac{1}{2}\Delta U(t,x)
+b(t,x)\cdot\nabla U(t,x)
-\lambda U(t,x)+b(t,x), \quad (t,x)\in (0,T]\times {\mathbb R}^n, \\
& U(0,x)=0, \quad  x\in{\mathbb R}^n,
\end{aligned}\right.
\end{align}
in the space
\begin{align*}
{\mathcal H}_T:=\Big\{\tilde U\in L^\infty([0,T];L^\infty_{\mathrm{loc}}(\mathbb{R}^n;\mathbb{R}^n)):\
|\nabla \tilde U|\in L^\infty([0,T]\times \mathbb{R}^n),\\
|\nabla^2 \tilde U|\in L^2([0,T];L^\infty(\mathbb{R}^n)), \ \
|\partial_t \tilde U|\in L^p([0,T];L^\infty_{\mathrm{loc}}(\mathbb{R}^n))\Big\},
\end{align*}
where $\lambda>0$ is chosen sufficiently large, and
$b\in L^p([0,T];{\mathcal C}^{\frac{2}{p}-1,\rho}_{l,d}(\mathbb{R}^n;\mathbb{R}^n))$ with $p\in (1,2]$.
A function $U$ is called a strong solution of \eqref{7.1} if
\[
U,\partial_t U,\partial^2_{x_i,x_j}U\in L^1([0,T];L^\infty_{\mathrm{loc}}(\mathbb{R}^n)),\qquad 1\le i,j\le n,
\]
it grows at most linearly in $x$, and \eqref{7.1} holds for a.e. $(t,x)\in[0,T]\times \mathbb{R}^n$.
We consider separately the cases $p\in(1,2)$ and $p=2$.

\medskip
\textbf{Case 1: $p\in (1,2)$.}
Extend $b$ from $[0,T]$ to $(-\infty,T]$ by setting $b(t,x)=0$ for $t<0$. Let $\varrho$ be a temporal mollifier:
\[
0\le \varrho\in {\mathcal C}^\infty_0(\mathbb{R}),\quad \mathrm{supp}(\varrho)\subset [0,1],\quad \int_{\mathbb{R}}\varrho(r)dr=1.
\]
For $m\in\mathbb{N}$, set $\varrho_m(r)=m\varrho(mr)$ and define the time-mollified field
\[
b_m(t,x)=(b(\cdot,x)\ast \varrho_m)(t)=\int_{\mathbb{R}} b(t-r,x)\varrho_m(r)dr.
\]
For $R>0$, let $\chi_R(x)=\chi(x/R)$ with
\begin{align}\label{7.2}
\chi\in {\mathcal C}^\infty_0(\mathbb{R}^n),\quad 0\le \chi\le 1,\quad |\nabla\chi|\le 2,\quad
\chi(x)=\begin{cases}
1,& {\text if} \ \  |x|<1,\\
0,& {\text if} \ \ |x|\geqq 2.
\end{cases}
\end{align}
Set $b_{m,R}(t,x)=b_m(t,x\chi_R(x))$. Then $b_{m,R}\in L^\infty([0,T];{\mathcal C}^{\frac{2}{p}-1,\rho}_{b,d}(\mathbb{R}^n;\mathbb{R}^n))$, and along a subsequence,
\begin{align}\label{7.3}
\lim_{m\to\infty}\;\lim_{R\to\infty} |b_{m,R}(t,x)-b(t,x)|=0,\quad \text{a.e. }(t,x)\in[0,T]\times\mathbb{R}^n.
\end{align}
Moreover,
\begin{align}\label{7.4}
\int_0^T \|(1+|\cdot|)^{-1} b_{m,R}(t,\cdot)\|_0^pdt\le \int_0^T \|(1+|\cdot|)^{-1} b(t,\cdot)\|_0^pdt
\end{align}
and
\begin{align}\label{7.5}
\begin{aligned}
{[b_{m,R}]}_{p,\frac{2}{p}-1,\rho}
&=\bigg(\int_0^T \sup_{0<|x-y|\le 1}\frac{|b_{m,R}(t,x)-b_{m,R}(t,y)|^p}{|x-y|^{2-p}\rho^p(|x-y|)}dt\bigg)^{\!\frac{1}{p}} \\
&\le 3^{\frac{2}{p}-1}[b]_{p,\frac{2}{p}-1,\rho}
\sup_{0<|x-y|\le 1}\frac{\rho(3|x-y|)}{\rho(|x-y|)}
\le C[b]_{p,\frac{2}{p}-1,\rho},
\end{aligned}
\end{align}
where
\[
[b_m]_{p,\frac{2}{p}-1,\rho}^p=\int_0^T[b_m(r,\cdot)]_{\frac{2}{p}-1,\rho}^p dr
\quad\text{and}\quad
[b]_{p,\frac{2}{p}-1,\rho}^p=\int_0^T[b(r,\cdot)]_{\frac{2}{p}-1,\rho}^pdr,
\]
and in \eqref{7.5} we used
\[
|x\chi_R(x)-y\chi_R(y)|\le 3|x-y|
\]
together with the slow variation of $\rho$ at $0$.

\smallskip
Consider the truncated Kolmogorov equation
\begin{align}\label{7.6}
\partial_t U_{m,R}(t,x)=\tfrac{1}{2}\Delta U_{m,R}(t,x) + b_{m,R}(t,x)\cdot\nabla U_{m,R}(t,x) - \lambda U_{m,R}(t,x) + b_{m,R}(t,x),
\end{align}
with $U_{m,R}(0,x)=0$. By Theorem~\ref{the3.6}\,\textup{(ii)}, there exists a unique
\[
U_{m,R}\in L^p([0,T];{\mathcal C}^{1+\frac{2}{p},\rho}_{l,d}(\mathbb{R}^n;\mathbb{R}^n))
\cap W^{1,p}([0,T];{\mathcal C}^{\frac{2}{p}-1,\rho}_{l,d}(\mathbb{R}^n;\mathbb{R}^n))
\cap L^\infty([0,T];{\mathcal C}^{1,\hat\rho}_{l,s}(\mathbb{R}^n;\mathbb{R}^n))
\]
solving \eqref{7.6}, where $\hat\rho$ is given by \eqref{3.10}. Moreover, the heat-kernel representation yields
\begin{align}\label{7.7}
U_{m,R}(t,x)=\int_0^t e^{-\lambda(t-r)}\big[K(t-r,\cdot)\ast \{b_{m,R}(r,\cdot)+b_{m,R}(r,\cdot)\cdot \nabla U_{m,R}(r,\cdot)\}\big](x)dr,
\end{align}
with $K(t,x)=(2\pi t)^{-n/2}e^{-|x|^2/(2t)}$.

\smallskip
Fix $x_0\in\mathbb{R}^n$ and consider
\begin{align}\label{7.8}
\dot x_t=-b_{m,R}(t,x_0+x_t),\qquad x_t|_{t=0}=0.
\end{align}
This ODE has a continuous solution. Define
\[
\hat U_{m,R}(t,x):=U_{m,R}(t,x+x_0+x_t),\qquad
\hat b_{m,R}(t,x):=b_{m,R}(t,x+x_0+x_t),
\]
and
\[
\tilde b_{m,R}(t,x):=\hat b_{m,R}(t,x)-\hat b_{m,R}(t,0).
\]
Then,
\begin{align}\label{7.9}
\begin{aligned}
\hat U_{m,R}(t,x)
&=\int_0^t e^{-\lambda(t-r)} \int_{\mathbb{R}^n} K(t-r,x-y)\tilde b_{m,R}(r,y)\cdot \nabla\hat U_{m,R}(r,y)dydr \\
&\quad + \int_0^t e^{-\lambda(t-r)}\int_{\mathbb{R}^n} K(t-r,x-y)\hat b_{m,R}(r,y)dydr.
\end{aligned}
\end{align}
Hence,
\begin{align}\label{7.10}
\begin{aligned}
|\nabla \hat U_{m,R}(t,0)|
&\le \int_0^t e^{-\lambda(t-r)} \int_{\mathbb{R}^n} |\nabla K(t-r,y)|[b_{m,R}(r,\cdot)]_{\frac{2}{p}-1,\rho} \\
&\qquad \times \big[1_{|y|\le 1}|y|^{\frac{2}{p}-1}\rho(|y|)+|y|1_{|y|>1}\big][1+\|\nabla U_{m,R}(r,\cdot)\|_0]dydr \\
&\le C\int_0^te^{-\lambda (t-r)}[b_{m,R}(r,\cdot)]_{\frac{2}{p}-1,\rho}[1+ \|\nabla U_{m,R}(r,\cdot)\|_0]\\ &\qquad \times\bigg(\int_0^1e^{-\frac{\tau^2}{4(t-r)}}(t-r)^{-\frac{n+1}{2}}\tau^{\frac{2}{p}+n-2}\rho(\tau)
d\tau+1\bigg)dr\\
&\le  C\int_0^te^{-\lambda (t-r)}[b_{m,R}(r,\cdot)]_{\frac{2}{p}-1,\rho}[1+ \|\nabla U_{m,R}(r,\cdot)\|_0]\\ &\qquad\times\bigg[\int_0^1e^{-\frac{\tau^2}{8(t-r)}}
\Big(\frac{\rho(\tau)}{\rho(\sqrt{t-r})}\Big)^{\frac{3p-2}{5p-2}}
\frac{\rho^{\frac{3p-2}{5p-2}}(\sqrt{t-r})}{(t-r)^{1-\frac{1}{p}}} \frac{\rho^{\frac{2p}{5p-2}}(\tau)}{\tau}
d\tau+1\bigg]dr,
\end{aligned}
\end{align}
where the last inequality results from the following estimate
\[
e^{-\frac{\tau^2}{4(t-r)}}(t-r)^{-\frac{n+1}{2}}\tau^{\frac{2}{p}+n-2}\rho(\tau)
\le
Ce^{-\frac{\tau^2}{8(t-r)}}\Big(\frac{\rho(\tau)}{\rho(\sqrt{t-r})}\Big)^{\frac{3p-2}{5p-2}}
\frac{\rho^{\frac{3p-2}{5p-2}}(\sqrt{t-r})}{(t-r)^{1-\frac{1}{p}}} \frac{\rho^{\frac{2p}{5p-2}}(\tau)}{\tau}.
\]

On the other hand,
\begin{align}\label{7.11}
\begin{aligned}
&\sup_{\tau\in (0,1], r\in (0,T]}\bigg[e^{-\frac{\tau^2}{8r}}
\Big(\frac{\rho(\tau)}{\rho(\sqrt{r})}
\Big)^{\frac{3p-2}{5p-2}}\bigg]  \\
&\le \sup_{r\in (0,T],\tau\in (0,\sqrt{r}] }\bigg[e^{-\frac{\tau^2}{8r}}
\Big(\frac{\rho(\tau)}{\rho(\sqrt{r})}
\Big)^{\frac{3p-2}{5p-2}}\bigg]
+\sup_{r\in (0,T],\tau\in[\sqrt{r},1]}\bigg[e^{-\frac{\tau^2}{8r}}
\Big(\frac{\rho(\tau)}{\rho(\sqrt{r})}
\Big)^{\frac{3p-2}{5p-2}}\bigg] \\
&\le 1+
\sup_{\mu\geqq 1}\sup_{r\in(0,1/\mu]}\bigg[e^{-\frac{\mu}{8}}\Big(\frac{\rho(\sqrt{\mu r})}{\rho(\sqrt{r})}
\Big)^{\frac{3p-2}{5p-2}}\bigg].
\end{aligned}
\end{align}
Choosing  $r_0=1$ in Lemma \ref{lem2.8}, in view of \eqref{2.7}, it follows that
\begin{align}\label{7.12}
\begin{aligned}
\sup_{r\in (0,1/\mu]}\frac{\rho(\sqrt{\mu r})}{\rho(\sqrt{r})}
 &= \sup_{r\in (0,1/\mu]}\exp\bigg\{c(\sqrt{\mu r})-c(\sqrt{r})+ \int_{\sqrt{r}}^{\sqrt{\mu r}}\frac{\zeta(\tau)}{\tau}d\tau\bigg\}\\
&\le \exp\Big\{2\sup_{0< \tau\le 1}|c(\tau)|+\sup_{0< \tau\le 1}\zeta(\tau)\log(\sqrt{\mu})\Big\}
\le C\mu^{\frac{1}{2}\sup\limits_{0< \tau\le 1}\zeta(\tau)}.
\end{aligned}
\end{align}
Combining \eqref{7.11} and \eqref{7.12}, one concludes
\[
\sup_{\tau\in [0,1], r\in (0,T]}\bigg[e^{-\frac{\tau^2}{8r}}
\Big(\frac{\rho(\tau)}{\rho(\sqrt{r})}
\Big)^{\frac{3p-2}{5p-2}}\bigg]
\le 1+C\sup_{\mu\geqq1}\Big[e^{-\frac{\mu}{8}}
\mu^{\frac{3p-2}{10p-4}\sup\limits_{0\le \tau\le 1}\zeta(\tau)}\Big]\le C.
\]

This, together with \eqref{7.10}, yields
\begin{align}\label{7.13}
|\nabla \hat{U}_{m,R}(t,0)|
\le C\int_0^t e^{-\lambda (t-r)}[b_{m,R}(r,\cdot)]_{\frac{2}{p}-1,\rho}\big[1+ \|\nabla U_{m,R}(r,\cdot)\|_0\big]\bigg[\frac{\rho^{\frac{3p-2}{5p-2}}(\sqrt{t-r})}{(t-r)^{1-\frac{1}{p}}}+1\bigg]dr,
\end{align}
where the last inequality uses the assumption that $\rho^{\frac{2p}{5p-2}}$ is a Dini function.

\smallskip
Since $x_0\in {\mathbb R}^n$ is arbitrary, we conclude from \eqref{7.13} that
\begin{align}\label{7.14}
\begin{aligned}
&\sup_{0\le t\le T}\|\nabla U_{m,R}(t,\cdot)\|_0\\  &\le
[b]_{p,\frac{2}{p}-1,\rho}[1+\sup_{0\le r\le T}\|\nabla U_{m,R}(r,\cdot)\|_0]
\bigg\{\int_0^{\sqrt{T}}e^{-\frac{\lambda pr^2}{p-1}}\Big[r+\frac{\rho^{\frac{p(3p-2)}{(5p-2)(p-1)}}(r)}{r}\Big]dr\bigg\}^{\frac{p-1}{p}}
\\ &\le C
[b]_{p,\frac{2}{p}-1,\rho}[1+\sup_{0\le r\le T}\|\nabla U_{m,R}(r,\cdot)\|_0]
\bigg\{\int_0^{\sqrt{T}}e^{-\frac{\lambda pr^2}{p-1}}\Big[r+\frac{\rho^{\frac{2p}{5p-2}}(r)}{r}\Big]dr\bigg\}^{\frac{p-1}{p}},
\end{aligned}
\end{align}
where the first inequality is due to the use of the H\"{o}lder inequality and \eqref{7.5}.

\smallskip
Since $\rho^{\frac{2p}{5p-2}}(r)/r\in L^1([0,\sqrt{T}])$, by dominated convergence, the integral in braces tends to $0$ as $\lambda\to\infty$. Choosing $\lambda$ sufficiently large so that
\[
C[b]_{p,\frac{2}{p}-1,\rho}
\bigg\{\int_0^{\sqrt{T}}e^{-\frac{\lambda pr^2}{p-1}}\Big[r+\frac{\rho^{\frac{2p}{5p-2}}(r)}{r}\Big]dr\bigg\}^{\frac{p-1}{p}}
\le \frac{1}{3},
\]
we obtain
\begin{align}\label{7.15}
\sup_{0\le t\le T}\|\nabla U_{m,R}(t,\cdot)\|_0
\le \frac{1}{2}.
\end{align}
This, together with \eqref{7.7} and \eqref{7.14}, also implies
\begin{align}\label{7.16}
\begin{aligned}
\sup_{0\le t\le T}\sup_{x\in {\mathbb R}^n}\frac{|U_{m,R}(t,x)|}{1+|x|}
&\le C\big[\|b_{m,1}\|_{L^p([0,T])}+\sup_{0\le t\le T}\|\nabla U_{m,R}(t,\cdot)\|_0\big]\\ &\le C\big[\|b_1\|_{L^p([0,T])}+[b]_{p,\frac{2}{p}-1,\rho}\big]
\le C\|b\|_{L^p([0,T];{\mathcal C}_{l,d}^{\frac{2}{p}-1,\rho}({\mathbb R}^n))},
\end{aligned}
\end{align}
where
\[
b_{m,1}(r)=\|(1+|\cdot|)^{-1}b_{m,R}(r,\cdot)\|_0,\qquad
b_1(r)=\|(1+|\cdot|)^{-1}b(r,\cdot)\|_0\in L^p([0,T]),
\]
and $\|b_{m,1}\|_{L^p([0,T])}\le \|b_1\|_{L^p([0,T])}$.

\smallskip
For the second-order derivatives of $U_{m,R}(t,x)$, we obtain an analogue of \eqref{7.13}:
\begin{align*}
|\nabla^2\hat{U}_{m,R}(t,0)|
&\le
C\int_0^te^{-\lambda (t-r)}[b_{m,R}(r,\cdot)]_{\frac{2}{p}-1,\rho}(t-r)^{-\frac{1}{2}}
\bigg[\frac{\rho^{\frac{3p-2}{5p-2}}(\sqrt{t-r})}{(t-r)^{1-\frac{1}{p}}}+1\bigg]dr,
\end{align*}
which, by Young's inequality, implies
\begin{align}\label{7.17}
\begin{aligned}
\|\nabla^2 U_{m,R}\|_{L^2([0,T];L^\infty({\mathbb R}^n))}
&\le C[b_{m,R}]_{p,\frac{2}{p}-1,\rho}\bigg\{\int_0^Te^{-\frac{2\lambda pr}{3p-2}}\Big[r^{-\frac{p}{3p-2}}+\frac{\rho^{\frac{2p}{5p-2}}(\sqrt{r})}{r}\Big]dr
\bigg\}^{\frac{3p-2}{2p}}\\
&\le C[b]_{p,\frac{2}{p}-1,\rho}\bigg\{1+\int_0^{\sqrt{T}}e^{-\frac{2\lambda pr^2}{3p-2}}\frac{\rho^{\frac{2p}{5p-2}}(r)}{r}dr\bigg\}^{\frac{3p-2}{2p}}
\le C[b]_{p,\frac{2}{p}-1,\rho}.
\end{aligned}
\end{align}

Combining \eqref{7.6}, \eqref{7.15}--\eqref{7.17} gives
\begin{align}\label{7.18}
\|(1+|\cdot|)^{-1}\partial_t U_{m,R}\|_{L^p([0,T];L^\infty(\mathbb{R}^n))}
\le C\|b\|_{L^p([0,T];{\mathcal C}^{\frac{2}{p}-1,\rho}_{l,d}(\mathbb{R}^n))}.
\end{align}

We now pass to the limit. Choose the diagonal approximation
$R=m$. For every $L>0$ and all sufficiently large $m$,
$b_{m,m}=b_m$ on $[0,T]\times\overline{B_L}$. Since
$b\in L^p([0,T];\mathcal C(\overline{B_L}))$, temporal mollification gives
\[
 b_{m,m}\longrightarrow b
 \quad\text{in }L^p([0,T];\mathcal C(\overline{B_L})), \quad {\text as} \quad m\rightarrow\infty.
\]

By \eqref{7.15}, \eqref{7.16} and \eqref{7.18},
\[
 |U_{m,m}(t,x)-U_{m,m}(s,y)|
 \le \frac12|x-y|+C_L|t-s|^{1-\frac1p},
 \qquad x,y\in\overline{B_L},\quad s,t\in[0,T].
\]
The Arzel\`a--Ascoli theorem and a diagonal subsequence
argument therefore yield a function $U$ such that
\begin{align}\label{7.19}
 U_{m,m}\longrightarrow U
 \quad\text{uniformly on }[0,T]\times\overline{B_L}
 \quad\text{for every }L>0, \quad {\text as} \quad m\rightarrow\infty.
\end{align}
In particular, $U(0,x)=0$.

\smallskip
To obtain convergence of the gradients, spatial interpolation
gives, for $w=U_{m,m}-U_{\ell,\ell}$,
\[
\begin{aligned}
 \|\nabla w\|_{L^4(0,T;L^\infty(B_L))}^4
\le C_L\Big(
 \|w\|_{L^\infty([0,T]\times B_{L+1})}^4
 +\|w\|_{L^\infty([0,T]\times B_{L+1})}^2
   \|\nabla^2w\|_{L^2([0,T];L^\infty(B_{L+1}))}^2
 \Big).
\end{aligned}
\]
Together with \eqref{7.17} and \eqref{7.19}, this shows that
$\nabla U_{m,m}$ is a Cauchy sequence in
$L^4([0,T];\mathcal C(\overline{B_L}))$. Its limit is the
distributional gradient of $U$. By \eqref{7.15} and
\eqref{7.19}, we obtain
\begin{align}\label{7.20}
 \nabla U_{m,m}\longrightarrow\nabla U
 \quad\text{in }L^q([0,T];\mathcal C(\overline{B_L})),
 \qquad 1\le q<\infty,  \quad {\text as} \quad m\rightarrow\infty.
\end{align}
After a further diagonal extraction, \eqref{7.17} and
\eqref{7.18} give
\begin{align}\label{7.21}
\begin{aligned}
 \nabla^2U_{m,m}\rightharpoonup\nabla^2U\ \text{in}\ L^2([0,T];L^2(B_L)),\
 \partial_tU_{m,m}\rightharpoonup\partial_tU\ \text{in}\quad L^p([0,T];L^2(B_L)), \ {\text as} \ m\rightarrow\infty.
\end{aligned}
\end{align}

Let $p'=p/(p-1)$. As $ m\rightarrow\infty$, the nonlinear term satisfies
\begin{align}\label{7.22}
\begin{aligned}
 &\|b_{m,m}\cdot\nabla U_{m,m}-b\cdot\nabla U\|_{L^1([0,T];\mathcal C(\overline{B_L}))}\\
 &\le
 C T^{1-\frac1p}
 \|b_{m,m}-b\|_{L^p([0,T];\mathcal C(\overline{B_L}))}+
 \|b\|_{L^p([0,T];\mathcal C(\overline{B_L}))}
 \|\nabla U_{m,m}-\nabla U\|
       _{L^{p'}(0,T;\mathcal C(\overline{B_L}))}
 \longrightarrow0.
\end{aligned}
\end{align}
Combining \eqref{7.19}--\eqref{7.22}, we can therefore pass
to the limit in \eqref{7.6} in the sense of distributions.
The resulting equation holds almost everywhere because
all its terms are locally integrable.

\smallskip
The bounds \eqref{7.16}--\eqref{7.18} pass to the limit by
lower semicontinuity on bounded spatial balls and then by
letting $L\to\infty$. Hence $U\in\mathcal H_T$.
The uniform spatial Lipschitz bound also passes to the limit,
giving
\begin{align}\label{7.23}
 \sup_{0\le t\le T}\|\nabla U(t,\cdot)\|_0
 \le\frac12.
\end{align}

Now we prove uniqueness. By linearity, it suffices to
consider the homogeneous problem, that is, with the source
term $b$ removed from \eqref{7.1}. Repeating the estimate
leading to \eqref{7.14} with $U$ in place of $U_{m,R}$ gives
\[
 \sup_{0\le t\le T}\|\nabla U(t,\cdot)\|_0
 \le \frac13
 \sup_{0\le t\le T}\|\nabla U(t,\cdot)\|_0.
\]
Hence $\nabla U\equiv0$. The homogeneous equation then
reduces to $\partial_tU=-\lambda U$ with $U(0,\cdot)=0$,
so $U\equiv0$.

\medskip
\textbf{Case 2: $p=2$.}
Let $0\le\widetilde\varrho\in
\mathcal C^\infty_0(\mathbb R^n)$ be a spatial mollifier with
$\mathrm{supp}(\widetilde\varrho)\subset B_1$ and
$\int_{\mathbb R^n}\widetilde\varrho(x)dx=1$.
For $k\in\mathbb N$, set
$\widetilde\varrho_k(x)=k^n\widetilde\varrho(kx)$ and define
\[
 b^k(t,x)=(b(t,\cdot)\ast\widetilde\varrho_k)(x),
 \qquad
 b^k_R(t,x)=b^k(t,x\chi_R(x)),
\]
with $\chi_R$ as in \eqref{7.2}. Then
\[
 b^k_R\in L^2([0,T];
 \mathcal C^{\beta,\rho}_{b,d}
 (\mathbb R^n;\mathbb R^n))
 \qquad\text{for every }\beta\in(0,1),
\]
and
\begin{align}\label{7.24}
 \lim_{k\to\infty}\;\lim_{R\to\infty}
 |b^k_R(t,x)-b(t,x)|=0,
 \quad\text{for a.e. }t\in[0,T]
 \text{ and every }x\in\mathbb R^n.
\end{align}
Moreover,
\begin{align}\label{7.25}
 \int_0^T
 \|(1+|\cdot|)^{-1}b^k_R(t,\cdot)\|_0^2dt
 \le
 4\int_0^T
 \|(1+|\cdot|)^{-1}b(t,\cdot)\|_0^2dt
\end{align}
and
\begin{align}\label{7.26}
 [b^k_R]_{2,0,\rho}\le C[b]_{2,0,\rho}.
\end{align}

Consider the Kolmogorov equation
\begin{align}\label{7.27}
\left\{
\begin{aligned}
 &\partial_tU^k_R(t,x)
 =\tfrac12\Delta U^k_R(t,x)
 +b^k_R(t,x)\cdot\nabla U^k_R(t,x)\\
 &\qquad\qquad\quad
 -\lambda U^k_R(t,x)+b^k_R(t,x),
 \quad (t,x)\in(0,T]\times\mathbb R^n,\\
 &U^k_R(0,x)=0,\quad x\in\mathbb R^n.
\end{aligned}
\right.
\end{align}
By Theorem~\ref{the3.6} \textup{(i)}, there exists a unique
\[
\begin{aligned}
 U^k_R\in{}&
 L^2([0,T];
 \mathcal C^{2+\beta,\rho}_{l,d}
 (\mathbb R^n;\mathbb R^n))\cap W^{1,2}([0,T];
 \mathcal C^{\beta,\rho}_{l,d}
 (\mathbb R^n;\mathbb R^n))\cap L^\infty([0,T];
 \mathcal C^{1+\beta,\rho}_{l,d}
 (\mathbb R^n;\mathbb R^n))
\end{aligned}
\]
solving \eqref{7.27}, and \eqref{7.7} holds with $b^k_R$
in place of $b_{m,R}$. Repeating the estimates for Case~1
yields
\begin{align}\label{7.28}
 \sup_{0\le t\le T}\|\nabla U^k_R(t,\cdot)\|_0
 \le\frac12
\end{align}
for sufficiently large $\lambda$, and
\begin{align}\label{7.29}
\left\{
\begin{aligned}
 &\sup_{0\le t\le T}\sup_{x\in\mathbb R^n}
 \frac{|U^k_R(t,x)|}{1+|x|}
 \le C\|b\|_{L^2([0,T];
 \mathcal C^{0,\rho}_{l,d}(\mathbb R^n))},\\[0.1cm]
 &\|\nabla^2U^k_R\|_{L^2([0,T];
 L^\infty(\mathbb R^n))}
 \le C[b]_{2,0,\rho},\\[0.1cm]
 &\|(1+|\cdot|)^{-1}\partial_tU^k_R\|
 _{L^2([0,T];L^\infty(\mathbb R^n))}
 \le C\|b\|_{L^2([0,T];
 \mathcal C^{0,\rho}_{l,d}(\mathbb R^n))}.
\end{aligned}
\right.
\end{align}

We now choose $R=k$. For every $L>0$ and all sufficiently
large $k$, $b^k_k=b^k$ on $[0,T]\times\overline{B_L}$.
Spatial mollification gives
\[
 \|b^k_k-b\|_{L^2([0,T];\mathcal C(\overline{B_L}))}
 \le C\rho(1/k)[b]_{2,0,\rho}
 \longrightarrow0
 \qquad\text{as }k\to\infty.
\]
The bounds \eqref{7.28}--\eqref{7.29} allow us to apply the
compactness argument following \eqref{7.18}, now with $p=2$.
After extracting a subsequence, we obtain
\[
\begin{aligned}
 U^k_k&\longrightarrow U\ \text{uniformly on}\ [0,T]\times\overline{B_L}\ \ \text{and} \ \
 \nabla U^k_k\longrightarrow\nabla U \ \text{in} \ L^q([0,T];{\mathcal C}(\overline{B_L})),
 \qquad 1\le q<\infty,
\end{aligned}
\]
together with weak convergence of the time derivatives
and Hessians in the corresponding local spaces.
In particular,
\[
 b^k_k\cdot\nabla U^k_k
 \longrightarrow b\cdot\nabla U
 \quad\text{in }L^1([0,T];\mathcal C(\overline{B_L})).
\]
Passing to the limit in \eqref{7.27} gives \eqref{7.1}.
Lower semicontinuity preserves \eqref{7.28}--\eqref{7.29},
so $U\in\mathcal H_T$ and \eqref{7.23} holds.
Uniqueness follows as in Case~1.

\medskip
\textbf{Part II:}
Well-posedness of the SDE \eqref{1.1} with low regularity,
linearly growing drift.

\smallskip
Apply Part~I with the time-reversed drift $b(T-t,x)$
in place of $b(t,x)$, and denote the corresponding solution
by $U$. Set $V(t,x)=U(T-t,x)$. Then $V\in\mathcal H_T$
satisfies
\begin{align}\label{7.30}
\left\{
\begin{aligned}
&\partial_t V(t,x)+\tfrac{1}{2}\Delta V(t,x)+b(t,x)\cdot\nabla V(t,x)
= \lambda V(t,x)-b(t,x),\quad (t,x)\in[0,T)\times\mathbb{R}^n,\\
&V(T,x)=0, \quad x\in \mathbb{R}^n.
\end{aligned}
\right.
\end{align}
Moreover, \eqref{7.23} holds for $V$.
Set $\Phi(t,x)=x+V(t,x)$. Then $\Phi(t,\cdot)$ is a
bi-Lipschitz homeomorphism uniformly in $t\in[0,T]$, and
\begin{align}\label{7.31}
\left\{
\begin{aligned}
&\frac{1}{2}\le \inf_{0\le t\le T,x\in \mathbb{R}^n}|\nabla\Phi(t,x)|\le \sup_{0\le t\le T,x\in \mathbb{R}^n}|\nabla\Phi(t,x)|\le \frac{3}{2},\\[0.1cm]
&\frac{2}{3}\le \inf_{0\le t\le T,x\in \mathbb{R}^n}|\nabla\Psi(t,x)|\le \sup_{0\le t\le T,x\in \mathbb{R}^n}|\nabla\Psi(t,x)|\le 2,
\end{aligned}
\right.
\end{align}
where $\Psi(t,\cdot)=\Phi(t,\cdot)^{-1}$.

\smallskip
For $0<\varepsilon<1$, define the time average
\[
 V_\varepsilon(t,x)
 =\frac1\varepsilon\int_t^{t+\varepsilon}V(r,x)dr
 =\int_0^1V(t+\varepsilon r,x)dr,
\]
and set
\[
 \Phi_\varepsilon(t,x)=x+V_\varepsilon(t,x),
 \qquad
 V(t,x)=V(T,x)=0\quad\text{for }t>T.
\]
Then
\[
 \Phi_\varepsilon\in
 W^{1,2}([0,T];
 W^{2,\infty}_{\mathrm{loc}}
 (\mathbb R^n;\mathbb R^n)).
\]
If $X_{s,t}(x)$ is a strong solution to \eqref{1.1},
It\^o's formula
(see, e.g., \cite[Theorem 3.7]{KR}) gives
\begin{align}\label{7.32}
\begin{aligned}
 \Phi_\varepsilon(t,X_{s,t}(x))
 &=\Phi_\varepsilon(s,x)
 +\int_s^t
 \partial_rV_\varepsilon(r,X_{s,r}(x))dr+
 \int_s^t
 b(r,X_{s,r}(x))\cdot
 \nabla V_\varepsilon(r,X_{s,r}(x))dr\\
 &\quad+
 \frac12\int_s^t
 \Delta V_\varepsilon(r,X_{s,r}(x))dr
 +\int_s^tb(r,X_{s,r}(x))dr\\
 &\quad+
 \int_s^t
 [I+\nabla V_\varepsilon(r,X_{s,r}(x))]dW_r.
\end{aligned}
\end{align}

We justify the passage to the limit in \eqref{7.32}.
Extend also $b$ by zero for $t>T$, and set
\[
\begin{aligned}
 b_\varepsilon(t,x)=\frac1\varepsilon\int_t^{t+\varepsilon}b(r,x)dr \quad \text{and} \quad
 (b\cdot\nabla V)_\varepsilon(t,x)=\frac1\varepsilon\int_t^{t+\varepsilon}
       b(r,x)\cdot\nabla V(r,x)dr.
\end{aligned}
\]
Since $V(T,x)=0$, the zero extension of $V$ introduces
no jump at $T$. Averaging \eqref{7.30} gives
\[
 \partial_tV_\varepsilon(t,x)
 +\frac12\Delta V_\varepsilon(t,x)
 =\lambda V_\varepsilon(t,x)-b_\varepsilon(t,x)
 -(b\cdot\nabla V)_\varepsilon(t,x).
\]
Consequently,
\begin{align}\label{7.33}
\begin{aligned}
 &\partial_tV_\varepsilon(t,x)
 +\frac12\Delta V_\varepsilon(t,x)
 +b(t,x)\cdot\nabla V_\varepsilon(t,x)
 +b(t,x)-\lambda V_\varepsilon(t,x)\\
 &=(b-b_\varepsilon)(t,x)
 +b(t,x)\cdot
   (\nabla V_\varepsilon(t,x)-\nabla V(t,x))+
 b(t,x)\cdot\nabla V(t,x)
 -(b\cdot\nabla V)_\varepsilon(t,x).
\end{aligned}
\end{align}
Fix $L>0$. The time regularity of $V$ implies
\begin{align}\label{7.34}
 \|V_\varepsilon-V\|
 _{L^\infty([0,T]\times B_L)}
 \longrightarrow0
 \qquad\text{as }\varepsilon\downarrow0.
\end{align}
For almost every $t$, $\nabla V(t,\cdot)$ has a Lipschitz
representative on $\overline{B_L}$. Since $\nabla V$ is
bounded, convergence of time averages gives
\begin{align}\label{7.35}
 \nabla V_\varepsilon\longrightarrow\nabla V
 \quad\text{in }L^q([0,T];\mathcal C(\overline{B_L})),
 \qquad 1\le q<\infty,
 \quad\text{as }\varepsilon\downarrow0.
\end{align}
Likewise,
\begin{align}\label{7.36}
b_\varepsilon\longrightarrow b,
\quad (b\cdot\nabla V)_\varepsilon\longrightarrow b\cdot\nabla V\quad
\text{in }L^p([0,T];{\mathcal C}(\overline{B_L})),\quad {\text as} \quad \varepsilon\downarrow 0.
\end{align}
Taking $p'=p/(p-1)$, H\"older's inequality yields
\begin{align}\label{7.37}
\begin{aligned}
 &\big\|
 (b-b_\varepsilon)
 +b\cdot(\nabla V_\varepsilon-\nabla V)
 +b\cdot\nabla V-(b\cdot\nabla V)_\varepsilon
 \big\|_{L^1([0,T];\mathcal C(\overline{B_L}))}\\
 &\le T^{1-\frac1p}
 \Big(
 \|b-b_\varepsilon\|_{L^p([0,T];\mathcal C(\overline{B_L}))}
 +\|b\cdot\nabla V-(b\cdot\nabla V)_\varepsilon\|
 _{L^p([0,T];\mathcal C(\overline{B_L}))}
 \Big)\\
 &\quad+
 \|b\|_{L^p([0,T];\mathcal C(\overline{B_L}))}
 \|\nabla V_\varepsilon-\nabla V\|
 _{L^{p'}(0,T;\mathcal C(\overline{B_L}))}
 \longrightarrow0
 \qquad\text{as }\varepsilon\downarrow0.
\end{aligned}
\end{align}

For $L>|x|$, define
\[
 \tau_L
 =\inf\{r\in[s,T]:|X_{s,r}(x)|\geq L\}\wedge T,
\]
with $\inf\varnothing=\infty$.
By the Burkholder--Davis--Gundy inequality,
\begin{align}\label{7.38}
\begin{aligned}
 &\mathbb E\sup_{s\le t\le T}
 \left|
 \int_s^{t\wedge\tau_L}
 (\nabla V_\varepsilon-\nabla V)
 (r,X_{s,r}(x))dW_r
 \right|^2\\
 &\le
 C\|\nabla V_\varepsilon-\nabla V\|
 _{L^2([0,T];\mathcal C(\overline{B_L}))}^2
 \longrightarrow0
 \qquad\text{as }\varepsilon\downarrow0.
\end{aligned}
\end{align}
On the other hand, the endpoint terms converge by local
uniform convergence of $V_\varepsilon$.
Combining \eqref{7.33}--\eqref{7.37} with \eqref{7.38},
we may therefore pass to the limit in \eqref{7.32},
stopped at $\tau_L$, with convergence in probability
uniformly in time. Finally, letting $L\to\infty$ gives,
almost surely for all $t\in[s,T]$,
\begin{align}\label{7.39}
\begin{aligned}
 \Phi(t,X_{s,t}(x))=\Phi(s,x)
 +\lambda\int_s^tV(r,X_{s,r}(x))dr+
 \int_s^t[I+\nabla V(r,X_{s,r}(x))]dW_r,
\end{aligned}
\end{align}
where the identity follows from \eqref{7.30}.

\smallskip
Denote $Y_{s,t}(y)=\Phi(t,X_{s,t}(x))$ with
$y=\Phi(s,x)$. Then
\begin{align}\label{7.40}
\left\{
\begin{aligned}
 &dY_{s,t}(y)
 =\widetilde b(t,Y_{s,t}(y))dt
 +\widetilde\sigma(t,Y_{s,t}(y))dW_t,
 \quad t\in(s,T],\\
 &Y_{s,s}(y)=y,
\end{aligned}
\right.
\end{align}
where
\[
 \widetilde b(t,z)
 =\lambda V\bigl(t,\Psi(t,z)\bigr),
 \qquad
 \widetilde\sigma(t,z)
 =I+\nabla V\bigl(t,\Psi(t,z)\bigr).
\]

Conversely, let $Y_{s,t}(y)$ solve \eqref{7.40}.
Differentiating $\Phi(t,\Psi(t,z))=z$ gives
\[
\begin{aligned}
& \ \ \ \nabla\Psi(t,z)
 =\bigl[I+\nabla V(t,\Psi(t,z))\bigr]^{-1},\quad \partial_t\Psi(t,z)=-\nabla\Psi(t,z)\,
       \partial_tV(t,\Psi(t,z)),\\
 &\nabla^2\Psi(t,z)[z_1,z_2]
 =-\nabla\Psi(t,z)\,
       \nabla^2V(t,\Psi(t,z))
       [\nabla\Psi(t,z)z_1,\nabla\Psi(t,z)z_2], \ z_1,z_2\in \mathbb{R}^n.
\end{aligned}
\]
The identities involving weak derivatives hold almost
everywhere. In particular,
\[
 \|\nabla\Psi(t,\cdot)\|_0\le2,
 \qquad
 \|\nabla^2\Psi(t,\cdot)\|_0
 \le8\|\nabla^2V(t,\cdot)\|_0,
\]
and $\partial_t\Psi\in
L^p([0,T];L^\infty_{\mathrm{loc}})$.
After localization, these bounds permit the generalized
It\^o formula, justified by temporal and spatial
regularization as in the argument leading to \eqref{7.39}.
Using \eqref{7.30}, we have
\[
\begin{aligned}
 \partial_t\Psi
 +\widetilde b\cdot\nabla\Psi
 +\frac12\sum_{i,j=1}^n
 (\widetilde\sigma\widetilde\sigma^\top)_{i,j}
       \partial^2_{z_i,z_j}\Psi
 =b(t,\Psi(t,z)),\ \ \nabla\Psi(t,z)\widetilde\sigma(t,z)=I.
\end{aligned}
\]
Consequently, $X_{s,t}(x)=\Psi(t,Y_{s,t}(y))$, with
$y=\Phi(s,x)$, solves \eqref{1.1}.
Hence \eqref{1.1} and \eqref{7.40} are equivalent.

\smallskip
The regularity of $V$ implies
\[
 \widetilde b\in
 L^\infty([0,T];
 \operatorname{Lip}(\mathbb R^n;\mathbb R^n))
 \quad\text{and}\quad
 \widetilde\sigma\in
 L^2([0,T];
 W^{1,\infty}(\mathbb R^n;\mathbb R^{n\times n})).
\]
Therefore, classical SDE theory yields a unique strong
solution $Y_{s,t}(y)$ to \eqref{7.40} together with the
flow property
\begin{equation}\label{7.41}
 Y_{s,t}(y)=Y_{r,t}(Y_{s,r}(y)),
 \qquad 0\le s\le r\le t\le T,
 \qquad Y_{s,s}(y)=y.
\end{equation}

Applying It\^o's formula to $|Y_{s,t}(y)|^q$ yields,
for $q\geq2$,
\[
\begin{aligned}
 d|Y_{s,t}(y)|^q
 \le C(1+|Y_{s,t}(y)|^q)dt+
 q|Y_{s,t}(y)|^{q-2}
 \langle Y_{s,t}(y),
 \widetilde\sigma(t,Y_{s,t}(y))dW_t\rangle,
\end{aligned}
\]
and the martingale term has zero expectation; hence
\begin{align}\label{7.42}
 \sup_{s\le t\le T}
 \mathbb E|Y_{s,t}(y)|^q
 \le C(1+|y|^q).
\end{align}

To prove that, almost surely, $Y_{s,t}(\cdot)$ is a
homeomorphism of $\mathbb R^n$ for every $0\le s\le t\le T$,
it is sufficient (cf. \cite[Lemmas II.2.4, II.4.1, II.4.2]{Kun84}) to show that for all $y,y'\in\mathbb{R}^n$ and $s<t$, $s'<t'$:
\begin{align}\label{7.43}
\sup_{s\le t\le T}\mathbb{E}|Y_{s,t}(y)-Y_{s,t}(y')|^{2\xi}\le C|y-y'|^{2\xi}\quad\text{for all }\xi<0,\quad y\ne y',
\end{align}
and
\begin{align}\label{7.44}
\mathbb{E}|Y_{s,t}(y)-Y_{s',t'}(y')|^q
\le C\Big\{|y-y'|^q + \big(1+|y|^q+|y'|^q\big)\big(|s-s'|^{\frac{q}{2}}+|t-t'|^{\frac{q}{2}}\big)\Big\},\quad q\geq 2.
\end{align}
For \eqref{7.43}, set
$f_\varepsilon(x)=\varepsilon+|x|^2$ and
$Y_{s,t}(y,y'):=Y_{s,t}(y)-Y_{s,t}(y')$.
It\^o's formula gives
\begin{align}\label{7.45}
\begin{aligned}
 f_\varepsilon^\xi(Y_{s,t}(y,y'))
 &\le f_\varepsilon^\xi(y-y')
 +C|\xi|\int_s^t
 f_\varepsilon^\xi(Y_{s,r}(y,y'))dr+
 C|\xi(\xi-1)|\int_s^t
 \kappa(r)^2f_\varepsilon^\xi(Y_{s,r}(y,y'))dr\\
 &\quad+
 2\xi\int_s^t
 f_\varepsilon^{\xi-1}(Y_{s,r}(y,y'))
 \big\langle Y_{s,r}(y,y'),
 \big(\widetilde\sigma(r,Y_{s,r}(y))
 -\widetilde\sigma(r,Y_{s,r}(y'))\big)dW_r
 \big\rangle,
\end{aligned}
\end{align}
with $\kappa(r)=\|\nabla^2V(r,\cdot)\|_0\in L^2([0,T])$.

\smallskip
For fixed $\varepsilon>0$, the integrand of the stochastic
integral is bounded in absolute value by
$C_\xi\varepsilon^\xi\kappa(r)$.
Since $\kappa\in L^2([0,T])$, this stochastic integral is a
square-integrable martingale and has zero expectation.
Taking expectations and applying Gr\"onwall's inequality yields
\[
\begin{aligned}
 \mathbb E
 \bigl(\varepsilon+|Y_{s,t}(y)-Y_{s,t}(y')|^2\bigr)^\xi\le
 (\varepsilon+|y-y'|^2)^\xi
 \exp\bigg(
 C_\xi\int_s^t(1+\kappa(r)^2)dr
 \bigg).
\end{aligned}
\]
The constant is independent of $\varepsilon$.
Since $\xi<0$, the integrand on the left increases
monotonically as $\varepsilon\downarrow0$.
For $y\ne y'$, the monotone convergence theorem therefore gives
\[
 \mathbb E|Y_{s,t}(y)-Y_{s,t}(y')|^{2\xi}
 \le C_\xi|y-y'|^{2\xi}.
\]
The constant is uniform for $s\le t\le T$, which proves
the required negative-moment estimate \eqref{7.43}.

\smallskip
To prove \eqref{7.44}, assume without loss of generality that $s<s'<t<t'$. Then
\begin{align}\label{7.46}
|Y_{s,t}(y)-Y_{s',t'}(y')|^q
\le 3^{q-1}\big\{|Y_{s,t}(y)-Y_{s,t}(y')|^q + |Y_{s,t}(y')-Y_{s',t}(y')|^q + |Y_{s',t}(y')-Y_{s',t'}(y')|^q\big\}.
\end{align}
The first term is controlled by
\begin{align}\label{7.47}
\sup_{s\le t\le T}\mathbb{E}|Y_{s,t}(y)-Y_{s,t}(y')|^q\le C|y-y'|^q.
\end{align}
For the second,
\begin{align}\label{7.48}
\mathbb{E}|Y_{s,t}(y')-Y_{s',t}(y')|^q
\le C(1+|y'|^q)|s-s'|^{\frac{q}{2}}.
\end{align}
For the third,
\begin{align}\label{7.49}
\mathbb{E}|Y_{s',t}(y')-Y_{s',t'}(y')|^q
\le C(1+|y'|^q)|t-t'|^{\frac{q}{2}}.
\end{align}
Combining \eqref{7.46}--\eqref{7.49} proves \eqref{7.44}. Hence $Y_{s,t}(\cdot)$ is a homeomorphism.

\smallskip
By \eqref{7.44} and Kolmogorov's continuity theorem, $Y$ has a
jointly continuous version. For this version, \eqref{7.41} holds
on a single event of probability one, first on a countable dense
set of admissible $(s,r,t,y)$ and then everywhere by continuity.
On this event, applying \eqref{7.41} twice gives
$Y_{0,T}=Y_{t,T}\circ Y_{s,t}\circ Y_{0,s}$. Consequently,
\begin{align}\label{7.50}
Y_{s,t}^{-1}(y)=Y_{0,s}\bigl(Y_{0,T}^{-1}(Y_{t,T}(y))\bigr).
\end{align}
Since $Y_{0,T}^{-1}$ is continuous in space and the forward maps
are jointly continuous, \eqref{7.50} gives the continuity of
$Y_{s,t}^{-1}(y)$ in $(s,t,y)$ a.s. Therefore, $\{Y_{s,t}(\cdot),\,t\in[s,T]\}$ forms a stochastic flow of homeomorphisms for \eqref{7.40}. By the equivalence between \eqref{1.1} and \eqref{7.40}, the family $\{X_{s,t}(\cdot),\,t\in[s,T]\}$ forms a stochastic flow of homeomorphisms solving \eqref{1.1}.

\smallskip
Let $R>0$ and let $x,x'\in\mathbb R^n$ satisfy
$|x|,|x'|\le R$.
For $0\le s\le t\le T$ and
$0\le s'\le t'\le T$,
\begin{align}\label{7.51}
\begin{aligned}
 |X_{s,t}(x)-X_{s',t'}(x')|
\le
 |\Psi(t,Y_{s,t}(y))-\Psi(t',Y_{s,t}(y))|+
 |\Psi(t',Y_{s,t}(y))
 -\Psi(t',Y_{s',t'}(y'))|,
\end{aligned}
\end{align}
where $y=\Phi(s,x)$ and $y'=\Phi(s',x')$.

\smallskip
Since $V$ satisfies equation \eqref{7.30}, we have
\[
\Phi\in W^{1,p}\big([0,T];L^\infty_{\mathrm{loc}}(\mathbb{R}^n;\mathbb{R}^n)\big),
\]
and $\Phi$ grows at most linearly in space. By \eqref{7.31} and the Sobolev embedding
\[
W^{1,p}([0,T])\hookrightarrow \mathcal{C}_b^{1-\frac{1}{p}}([0,T]),
\]
we obtain
\begin{align*}
\begin{aligned}
|\Psi(t,Y_{s,t}(y))-\Psi(t^\prime,Y_{s,t}(y))|
&\le 2\big|\Phi(t,\Psi(t,Y_{s,t}(y)))-\Phi(t,\Psi(t^\prime,Y_{s,t}(y)))\big| \\
&= 2\big|\Phi(t^\prime,\Psi(t^\prime,Y_{s,t}(y)))-\Phi(t,\Psi(t^\prime,Y_{s,t}(y)))\big| \\
&\le C|t-t^\prime|^{1-\frac{1}{p}}(1+|Y_{s,t}(y)|).
\end{aligned}
\end{align*}
This, together with \eqref{7.51}, \eqref{7.31},
\eqref{7.42} and \eqref{7.44}, yields, for every $q\geq2$,
\begin{align}\label{7.52}
\begin{aligned}
 \mathbb E|X_{s,t}(x)-X_{s',t'}(x')|^q
 &\le
 C|t-t'|^{(1-\frac1p)q}
 \bigl(1+\mathbb E|Y_{s,t}(y)|^q\bigr)+
 C\mathbb E|Y_{s,t}(y)-Y_{s',t'}(y')|^q\\
 &\le C\Big(
 |t-t'|^{(1-\frac1p)q}
 +|y-y'|^q
 +|s-s'|^{q/2}+|t-t'|^{q/2}
 \Big).
\end{aligned}
\end{align}
Notice that
\[
\begin{aligned}
 |y-y'|\le
 |\Phi(s,x)-\Phi(s',x)|
 +|\Phi(s',x)-\Phi(s',x')|\le
 C\bigl(|s-s'|^{1-\frac1p}+|x-x'|\bigr).
\end{aligned}
\]
Since $p\in(1,2]$, we conclude that
\[
 \mathbb E|X_{s,t}(x)-X_{s',t'}(x')|^q
 \le C\bigl(
 |x-x'|^q
 +|s-s'|^{(1-\frac1p)q}
 +|t-t'|^{(1-\frac1p)q}
 \bigr).
\]
This completes the proof. \qed

\section{Conclusions}\label{sec8}\setcounter{equation}{0}
In this paper, we studied the stochastic differential equation (SDE) \eqref{1.1} with critically low regularity, linearly growing drift in the Lebesgue--H\"older--Dini framework. Our main result shows that, under the assumptions
\[
b\in L^p([0,T];{\mathcal C}^{\frac{2}{p}-1,\rho}_{l,d}(\mathbb{R}^n;\mathbb{R}^n)),
\qquad p\in(1,2],
\]
where $\rho$ is slowly varying at zero and $\rho^{\frac{2p}{5p-2}}$ is a Dini function, the SDE \eqref{1.1} admits a unique stochastic flow of homeomorphisms. In particular, in the borderline case $p=2$, this result provides a partial answer to the open problem concerning strong well-posedness and flow generation for SDEs with drifts in a critical $L^2_tL^\infty_x$-type regime.

\smallskip
The proof relies on a maximal regularity theory for parabolic equations in Lebesgue--H\"older--Dini spaces developed in this paper. By combining refined heat-kernel estimates with a vector-valued Calder\'on--Zygmund theorem, we obtained mixed space-time regularity estimates for parabolic equations with low regularity coefficients. These analytic estimates serve as the key input for the It\^o--Tanaka trick (or Zvonkin transform), which allows us to reduce the original SDE to an auxiliary equation with Lipschitz coefficients.

\smallskip
It is worth noting that the maximal regularity results established in this work apply to a broader class of parabolic equations, including those with non-constant diffusion coefficients. However, in the stochastic part of the paper we restrict ourselves to the case of constant diffusion for the SDE \eqref{1.1}. Extending the stochastic flow results to SDEs with variable or degenerate diffusion coefficients under the same critical Lebesgue--H\"older--Dini framework remains an interesting and challenging direction for future research.

\smallskip
Overall, the results of this paper suggest that the Lebesgue--H\"older--Dini scale provides a natural borderline framework for studying stochastic flows of SDEs with rough drifts, and that maximal regularity techniques offer an effective analytic approach for understanding stochastic dynamics under minimal spatial regularity assumptions.

\section*{Acknowledgements} This work was supported in part by National Natural Science Foundation of China grants 12371243, 11771123, 12171247, the Startup Foundation for Introducing Talent of NUIST, Jiangsu Provincial Double-Innovation Doctor Program JSSCBS20210466 and Qing Lan Project.

\bigskip\noindent
\textbf{Conflicts of Interest} The authors declare that they have no competing interests.

\bigskip\noindent
\textbf{Data Availability Statements} Data sharing is not applicable to this article as no new data were created or analyzed in this study.

\end{document}